\input amstex
\input xy
\xyoption{all}
\documentstyle{amsppt}
\document

\magnification 1000

\def\gen{\frak{g}}

\def\ben{\frak{b}}

\def\len{\frak{l}}

\def\Sen{\frak{S}}

\def\g{{\gamma}}

\def\l{{{\lambda}}}

\def\eps{{\varepsilon}}

\def\Bcb{{\Bb}}

\def\eb{{\bold e}}
\def\fb{{\bold f}}

\def\ib{{\bold i}}

\def\jb{{\bold j}}
\def\kb{{\bold k}}

\def\forb{{{\bold{for}}}}

\def\Ab{{\bold A}}

\def\Bb{{\bold B}}

\def\Fb{{\bold F}}

\def\Modb{{\bold{Mod}}}
\def\fModb{{\bold{fMod}}}
\def\modb{{\bold{mod}}}
\def\fmodb{{\bold{fmod}}}

\def\projb{{\bold{proj}}}

\def\Gb{{\bold G}}

\def\Kb{{\bold K}}

\def\Nb{{\bold N}}
\def\Hb{{\bold H}}

\def\Rb{{\bold R}}

\def\Sb{{\bold S}}

\def\Lb{{\bold L}}

\def\Vb{{\bold V}}

\def\k{{\kb}}

\def\B{{\roman B}}

\def\A{{\roman A}}
\def\B{{\roman B}}
\def\D{{\roman D}}

\def\th{{{}^\theta}}
\def\cc{{{}^\circ}}

\def\ad{{\roman{ad}}}

\def\res{{\roman{res}}}
\def\top{{\text{top}}}
\def\gr{\roman{gr}}

\def\Ext{\roman{Ext}}

\def\Hom{{\roman{Hom}}}
\def\hom{{\roman{hom}}}
\def\gHom{{\roman{hom}}}

\def\Res{\roman{Res}}

\def\dim{{\roman{dim}}}
\def\gdim{{\roman{gdim}}}

\def\Ind{{\roman{Ind}}}

\def\End{{\roman{End}}}
\def\Hom{{\roman{Hom}}}

\def\Id{{\roman{id}}}

\def\proj{{{{{\projb}}}}}

\def\NN{{\Bbb N}}

\def\QQ{{\Bbb Q}}

\def\ZZ{{\Bbb Z}}

\def\Ac{{\Cal A}}

\def\Cc{{\Cal C}}

\def\Kc{{\Cal K}}
\def\Lc{{\Cal L}}

\def\Rc{{\Cal R}}

\def\Vc{{\Cal V}}
\def\Vcb{{\pmb{\Vc}}}

\def\Ccb{{\pmb{\Cc}}}

\def\mod{{{\roman{mod}}}}

\def\soc{\roman{soc}}
\def\top{\roman{top}}

\def\and{{\text{and}}}
\def\low{{\text{low}}}
\def\up{{\text{up}}}

\def\res{{\roman{res}}}

\def\max{{\roman{max}}}

\def\ds{\displaystyle}

\def\qed{\hfill $\sqcap \hskip-6.5pt \sqcup$}        % White box
\overfullrule=0pt                                    % No black boxes

\def\7dag{{{\!\!\!\!\!\!\!\dag}}}
\def\6dag{{{\!\!\!\!\!\!\dag}}}
\def\5dag{{{\!\!\!\!\!\dag}}}
\def\4dag{{{\!\!\!\!\dag}}}
\def\3dag{{{\!\!\!\dag}}}
\def\2dag{{{\!\!\dag}}}
\def\1dag{{{\!\dag}}}

\def\la{{\langle}}
\def\ra{{\rangle}}

\newdimen\Squaresize\Squaresize=14pt
\newdimen\Thickness\Thickness=0.5pt
\def\Square#1{\hbox{\vrule width\Thickness
          \alphaox to \Squaresize{\hrule height \Thickness\vss
          \hbox to \Squaresize{\hss#1\hss}
          \vss\hrule height\Thickness}
          \unskip\vrule width \Thickness}
          \kern-\Thickness}
\def\Vsquare#1{\alphaox{\Square{$#1$}}\kern-\Thickness}

\nologo

%%%%%%%%%%%%%%%%%%%%%%%%%%%%%%%%%%%%%%%%%%%%%%%%%%%%%%%%%%%%%%%%%%%%%%%%%%%%%%%%
\def\upd{{\cc}}

\def\leqs{\leqslant}
\def\geqs{\geqslant}
\def\indu{\roman{ind}}

%%%%%%%%%%%%%%%%%%%%%%%%%%%%%%%%%%%%%%%%%%%%%%%%%%%%%%%%%%%%%%%%%%%%%%%%%%%%%%%%
\topmatter
\title Canonical bases and affine Hecke algebras of type D
\endtitle
\rightheadtext{} \abstract
We prove a conjecture of Miemietz and Kashiwara
on canonical bases and branching rules of affine
Hecke algebras of type D. The proof is similar
to the proof of the type B case in \cite{VV}.
\endabstract
\author P. Shan, M. Varagnolo, E. Vasserot\endauthor
\address D\'epartement de Math\'ematiques,
Universit\'e Paris 7, 175 rue du Chevaleret, F-75013 Paris,
Fax : 01 44 27 78 18
\endaddress
\email shan\@math.jussieu.fr\endemail
\address D\'epartement de Math\'ematiques,
Universit\'e de Cergy-Pontoise, UMR CNRS 8088, F-95000 
Cergy-Pontoise,  Fax : 01 34 25 66 45\endaddress \email
michela.varagnolo\@math.u-cergy.fr\endemail
\address D\'epartement de Math\'ematiques,
Universit\'e Paris 7, 175 rue du Chevaleret, F-75013 Paris,
Fax : 01 44 27 78 18
\endaddress
\email vasserot\@math.jussieu.fr\endemail
\thanks
2000{\it Mathematics Subject Classification.} Primary ??; Secondary
??.
\endthanks
\endtopmatter
\document

\head Introduction\endhead
Let $\fb$ be the negative part of the quantized enveloping algebra of type $\A^{(1)}$.
Lusztig's description of the canonical basis of
$\fb$  implies that this basis can be naturally 
identified with the set of isomorphism classes of 
simple objects of a category of modules of the affine Hecke algebras of type A.
This identification was mentioned in \cite{G}, and  was used
in \cite{A}. 
More precisely, there is a linear isomorphism between 
$\fb$ and the
Grothendieck group of finite dimensional modules of the
affine Hecke algebras of type A, and it is proved in \cite{A}
that the induction/restriction functors
for affine Hecke algebras are given by the action of the Chevalley generators  and 
their transposed operators with respect to 
some symmetric bilinear form on $\fb$. 
 
The branching rules
for affine Hecke algebras of type B have been investigated
quite recently, see \cite{E}, \cite{EK1,2,3}, \cite{M} and \cite{VV}. 
In particular, in \cite{E}, \cite{EK1,2,3} an analogue of Ariki's construction is conjectured and 
studied for affine Hecke algebras of type B. Here $\fb$ is replaced by
a module $\th\Vb(\l)$ over an algebra $\th\!\Bcb$.
More precisely it is
conjectured there that $\th\Vb(\l)$ admits a canonical basis which is naturally identified with
the set of isomorphism classes of 
simple objects of a category of modules of the affine Hecke algebras of type B.
Further, in this identification the 
branching rules of the affine Hecke algebras of type B
should be given by the $\th\Bcb$-action on
$\th\Vb(\l)$. This conjecture has been proved \cite{VV}.
It uses both the geometric picture introduced in \cite{E} (to prove part of the conjecture)
and a new kind of graded algebras similar to the KLR algebras from \cite{KL}, \cite{R}.

A similar description of the branching rules
for affine Hecke algebras of type D has also been conjectured in \cite{KM}.
In this case $\fb$ is replaced by
another module $\cc\Vb$ over the algebra $\th\!\Bcb$ (the same algebra as in the type B case).
The purpose of this paper is to prove the type D case.
The method of the proof is the same as in \cite{VV}.
First we introduce a family of graded algebras $\cc\Rb_m$ for $m$ a non negative integer.
They can be viewed as the Ext-algebras of some  complex
of constructible sheaves naturally attached to the Lie algebra of the group $SO(2m)$,
see Remark 2.8.
This complex enters in the Kazhdan-Lusztig classification of the simple modules
of the affine Hecke algebra of the group $Spin(2m)$.
Then we identify $\cc\Vb$ with the sum of the
Grothendieck groups of the graded algebras $\cc\Rb_m$.

The plan of the paper is the following.
In Section 1 we  introduce  a graded algebra $\cc\Rb(\Gamma)_\nu$.
It is associated with a quiver $\Gamma$ with an involution $\theta$ and with a dimension vector $\nu$.
In Section 2 we consider a particular choice of pair $(\Gamma,\theta)$.
The graded algebras $\cc\Rb(\Gamma)_{\nu}$ associated with this pair $(\Gamma,\theta)$
are denoted by the symbol $\cc\Rb_m$.
Next we introduce the affine Hecke algebra of type D, more precisely the affine Hecke algebra associated
with the group $SO(2m)$, 
and we prove that it is Morita equivalent
to $\cc\Rb_m$.
In Section 3 we categorify the module $\cc\Vb$ from \cite{KM} using the graded algebras $\cc\Rb_m$,
see Theorem 3.28.
The main result of the paper is Theorem 3.33.

\head 0.~ Notation\endhead

\subhead 0.1.~Graded modules over graded algebras\endsubhead Let
$\kb$ be an algebraically closed
field of characteristic 0. By a graded $\kb$-algebra
$\Rb=\bigoplus_d\Rb_d$ we'll always mean a $\ZZ$-graded associative
$\kb$-algebra. Let $\Rb$-$\modb$ be the category of finitely
generated graded $\Rb$-modules, $\Rb$-$\fmodb$ be the full
subcategory of finite-dimensional graded modules and $\Rb$-$\proj$
be the full subcategory of projective objects. Unless specified
otherwise all modules are left modules. We'll abbreviate
$$K(\Rb)=[\Rb\text{-}\proj],\quad
G(\Rb)=[\Rb\text{-}\fmodb].$$ Here $[\Ccb]$ denotes the
Grothendieck group of an exact category $\Ccb$. Assume that the
$\kb$-vector spaces $\Rb_d$ are finite dimensional for each $d$.
Then $K(\Rb)$ is a free Abelian group with a basis formed by the
isomorphism classes of the indecomposable objects in $\Rb$-$\proj$,
and $G(\Rb)$ is a free Abelian group with a basis formed by the
isomorphism classes of the simple objects in $\Rb$-$\fmodb$. Given
an object $M$ of $\Rb$-$\proj$ or $\Rb$-$\fmodb$ let $[M]$ denote
its class in $K(\Rb)$, $G(\Rb)$ respectively. When there is no risk
of confusion we abbreviate $M=[M]$. We'll write $[M:N]$ for the
composition multiplicity of the $\Rb$-module $N$ in the $\Rb$-module
$M$. Consider the ring $\Ac=\ZZ[v,v^{-1}]$. If the grading of $\Rb$
is bounded below then the $\Ac$-modules $K(\Rb)$, $G(\Rb)$ are free.
Here $\Ac$ acts on $G(\Rb)$, $K(\Rb)$ as follows
$$vM=M[1],\quad v^{-1}M=M[-1].$$
For any $M,N$ in $\Rb$-$\modb$ let
$$\gHom_\Rb(M,N)=\bigoplus_d\Hom_\Rb(M,N[d])$$
be the $\ZZ$-graded $\kb$-vector space of all $\Rb$-module
homomorphisms. If $\Rb=\kb$ we'll omit the subscript $\Rb$ in hom's
and in tensor products. For any graded $\kb$-vector space $M=\bigoplus_dM_d$ we'll write
$$\gdim(M)=\sum_dv^d\dim(M_d),$$
where $\dim$ is the dimension over $\kb$.

\vskip3mm

\subhead 0.2.~Quivers with involutions\endsubhead Recall that a
quiver $\Gamma$ is a tuple $(I,H,h\mapsto h',h\mapsto h'')$ where
$I$ is the set of vertices, $H$ is the set of arrows and for each
$h\in H$ the vertices $h',h''\in I$ are the origin and the goal of
$h$ respectively. Note that the set $I$ may be infinite. We'll
assume that no arrow may join a vertex to itself. For each $i,j\in
I$ we write
$$H_{i,j}=\{h\in H;h'=i,h''=j\}.$$ We'll abbreviate $i\to j$ if
$H_{i,j}\neq\emptyset$. Let $h_{i,j}$ be the number of elements in
$H_{i,j}$ and set
$$i\cdot j=-h_{i,j}-h_{j,i},\quad i\cdot i=2,\quad i\neq j.$$
An involution $\theta$ on $\Gamma$ is a pair of involutions on $I$
and $H$, both denoted by $\theta$, such that
%$\theta(h)'=\theta(h'')$ and $\theta(h)''=\theta(h')$
the following properties hold for each $h$ in $H$ \vskip2mm
\itemitem{$\bullet$}
$\theta(h)'=\theta(h'')$ and $\theta(h)''=\theta(h')$, \vskip2mm
\itemitem{$\bullet$}
$\theta(h')=h''$ iff $\theta(h)=h$. \vskip2mm \noindent {\it We'll
always assume that $\theta$ has no fixed points in $I$}, i.e., there
is no $i\in I$ such that $\theta(i)=i$. To simplify we'll say that
$\theta$ has no fixed point. Let $$\th\NN I=\{\nu=\sum_i\nu_ii\in\NN
I: \nu_{\theta(i)}=\nu_i,\ \forall i\}.$$ For any $\nu\in\th\NN I$
set $|\nu|=\sum_i\nu_i$. It is an even integer. Write $|\nu|=2m$
with $m\in\NN$. We'll denote by $\th\! I^\nu$ the set of sequences
$$\ib=(i_{1-m},\ldots,i_{m-1},i_m)$$
of elements in $I$ such that $\theta(i_l)=i_{1-l}$ and
$\sum_ki_k=\nu$. For any such sequence $\ib$ we'll abbreviate
$\theta(\ib)=(\theta(i_{1-m}),\ldots,\theta(i_{m-1}),\theta(i_m))$.
Finally, we set
$$\th\! I^m=\bigcup_\nu\th\! I^\nu,\quad\nu\in\th\NN I,\quad |\nu|=2m.$$

\subhead 0.3.~The wreath product\endsubhead Given a positive integer
$m$, let $\Sen_m$ be the symmetric group, and $\ZZ_2=\{-1,1\}$.
Consider the wreath product $W_m=\Sen_m\wr\ZZ_2.$ Write
$s_1,\ldots,s_{m-1}$ for the simple reflections in $\Sen_m$. For
each $l=1,2,\dots m$ let $\eps_l\in(\ZZ_2)^m$ be $-1$ placed at the
$l$-th position. There is a unique action of $W_m$ on the set
$\{1-m,\dots,m-1,m\}$ such that $\Sen_m$ permutes $1,2,\dots m$ and
such that $\eps_l$ fixes $k$ if $k\neq l, 1-l$ and switches $l$ and
$1-l$. The group $W_m$ acts also on $\th\!I^\nu$. Indeed, view a
sequence $\ib$ as the map
$$\{1-m,\dots, m-1,m\}\to I,\quad l\mapsto i_l.$$
Then we set $w(\ib)=\ib\circ w^{-1}$ for $w\in W_m$. For each $\nu$ we fix once for all a sequence
$$\ib_e=(i_{1-m},\ldots,i_m)\in\th\! I^\nu.$$ Let $W_e$ be the
centralizer of $\ib_e$ in $W_m$. Then there is a bijection
$$W_e\backslash W_m\to\th\! I^\nu,\quad W_e w\mapsto
w^{-1}(\ib_e).$$
Now, assume that $m>1$. We set $s_0=\eps_1s_1\eps_1$. Let
$\cc W_m$ be the subgroup of $W_m$ generated by
$s_0,\ldots,s_{m-1}$. We'll regard it as a Weyl group of type $\D_m$
such that $s_0,\ldots,s_{m-1}$ are the simple reflections. Note that
$W_e$ is a subgroup of $\cc W_m$. Indeed, if $W_e\not\subset\cc W_m$ there should exist $l$ such that
$\eps_l$ belongs to $W_e$. This would imply that
$i_l=\theta(i_l)$, contradicting the fact that
$\theta$ has no fixed point. Therefore $\th\!I^\nu$ decomposes into two $\cc W_m$-orbits. We'll
denote them by $\th\!I^\nu_+$ and $\th\!I^\nu_-$. For $m=1$ we set
$\cc W_1=\{ e\} $ and we choose again $\th\!I^\nu_+$ and $\th\!I^\nu_-$ in a obvious way.

\vskip3mm

\head 1.~The graded $\kb$-algebra $\upd\Rb(\Gamma)_{\nu}$\endhead

Fix a quiver $\Gamma$ with set of vertices $I$ and set of arrows
$H$. Fix an involution $\theta$ on $\Gamma$. Assume that $\Gamma$
has no 1-loops and that $\theta$ has no fixed points. Fix a
dimension vector $\nu\neq 0$ in $\th\NN I$. Set $|\nu|=2m$.

\vskip3mm

\subhead 1.1.~Definition of the graded $\kb$-algebra
$\upd\Rb(\Gamma)_\nu$\endsubhead Assume that $m>1$. We define a
graded $\kb$-algebra $\upd\Rb(\Gamma)_\nu$ with $1$ generated by
$1_\ib$, $\varkappa_{l}$, $\sigma_{k}$, with $\ib\in\th\! I^\nu$,
$l=1,2,\dots,m$, $k=0,1,\ldots, m-1$ modulo the following defining
relations

\vskip2mm

\itemitem{$(a)$}
$1_\ib\,1_{\ib'}=\delta_{\ib,\ib'}1_\ib$, \quad $\sigma_{k}1_\ib=
1_{s_k\ib}\sigma_{k}$, \quad $\varkappa_{l}1_\ib=
1_{\ib}\varkappa_{l}$,

\vskip2mm

\itemitem{$(b)$}
$\varkappa_{l}\varkappa_{l'}=\varkappa_{l'}\varkappa_{l}$,

\vskip2mm

\itemitem{$(c)$}
$\sigma_{k}^21_\ib=
Q_{i_k,i_{s_k(k)}}(\varkappa_{s_k(k)},\varkappa_{k})1_\ib$,

\vskip2mm

\itemitem{$(d)$}
$\sigma_{k}\sigma_{k'}=\sigma_{k'}\sigma_{k}$ if $1\leqs k<k'-1<m-1$
or $0=k<k'\neq 2$,

\vskip2mm

\itemitem{$(e)$}
$(\sigma_{s_k(k)}\sigma_{k}\sigma_{s_k(k)}-
\sigma_{k}\sigma_{s_k(k)}\sigma_{k})1_\ib=$
$$=\cases
{\ds{Q_{i_k,i_{s_k(k)}}(\varkappa_{s_k(k)},\varkappa_{k})
-Q_{i_k,i_{s_k(k)}}(\varkappa_{s_k(k)},\varkappa_{s_k(k)+1})\over
\varkappa_{k}-\varkappa_{s_k(k)+1}}1_\ib}& \roman{if}\
i_k=i_{s_k(k)+1}, \vspace{2mm} 0&\roman{else},\endcases$$

\vskip2mm

\itemitem{$(f)$}
$(\sigma_{k}\varkappa_{l}-\varkappa_{s_k(l)}\sigma_{k})1_\ib =\cases
-1_\ib&\ \roman{if}\  l=k,\, i_k=i_{s_k(k)}, \cr 1_\ib&\ \roman{if}\
l=s_k(k),\, i_k=i_{s_k(k)},\cr 0&\ \roman{else}.
\endcases$

\vskip2mm

\noindent Here we have set $\varkappa_{1-l}=-\varkappa_l$ and
$$Q_{i,j}(u,v)=\cases(-1)^{h_{i,j}}(u-v)^{-i\cdot j}&\ \roman{if}\
i\neq j, \cr 0&\ \roman{else}.\endcases\leqno(1.1)$$
If $m=0$ we set $\cc\Rb(\Gamma)_0=\kb\oplus\kb$.
If $m=1$ then we have $\nu=i+\theta(i)$ for some $i\in I$. Write
$\ib=i\theta(i)$, and
$$\cc\Rb(\Gamma)_\nu=
\kb[\varkappa_1]1_{\ib}\oplus\kb[\varkappa_1]1_{\theta(\ib)}.$$
We'll abbreviate $\sigma_{\ib,k}=\sigma_{k}1_\ib$ and 
$\varkappa_{\ib,l}=\varkappa_{l}1_\ib$. The grading on
$\cc\Rb(\Gamma)_{0}$ is the trivial one. For $m\geqs 1$ the grading
on $\cc\Rb(\Gamma)_{\nu}$ is given by the following rules :
$$\aligned
&\deg(1_{\ib})=0,\vspace{2mm}
&\deg(\varkappa_{\ib,l})=2,\vspace{2mm}
&\deg(\sigma_{\ib,k})=-i_k\cdot i_{s_k(k)}.
\endaligned$$
We define $\omega$ to be the unique involution of the graded
$\kb$-algebra $\cc\Rb(\Gamma)_\nu$ which fixes $1_\ib$,
$\varkappa_l$, $\sigma_k$. We set $\omega$ to be identity on
$\cc\Rb(\Gamma)_0$.

\vskip3mm

\subhead 1.2.~Relation with the graded $\kb$-algebra
$\th\Rb(\Gamma)_\nu$\endsubhead A family of graded $\kb$-algebra
$\th\Rb(\Gamma)_{\l,\nu}$ was introduced in \cite{VV, sec.~ 5.1}, for $\l$
an arbitrary dimension vector in $\NN I$. Here we'll only consider the special
case $\l=0$, and we abbreviate
$\th\Rb(\Gamma)_\nu=\th\Rb(\Gamma)_{0,\nu}$.
Recall that if $\nu\neq 0$ then $\th\Rb(\Gamma)_\nu$ is the graded $\kb$-algebra with
$1$ generated by $1_\ib$, $\varkappa_{l}$, $\sigma_{k}$, $\pi_1$,
with $\ib\in\th\! I^\nu$, $l=1,2,\dots,m$, $k=1,\ldots, m-1$ such
that $1_\ib$, $\varkappa_{l}$ and $\sigma_{k}$ satisfy the same
relations as before and
$$\gathered \pi_1^2=1, \quad\pi_11_\ib\pi_1=1_{\eps_1\ib},\quad
\pi_1\varkappa_l\pi_1=\varkappa_{\eps_1(l)},\vspace{2mm}
(\pi_1\sigma_1)^2=(\sigma_1\pi_1)^2,\quad\pi_1\sigma_k\pi_1=\sigma_k\text{
if }k\neq 1.\endgathered$$ If $\nu=0$ then $\th\Rb(\Gamma)_0=\kb$. 
The grading is given by setting $\deg(1_{\ib})$,
$\deg(\varkappa_{\ib,l})$, $\deg(\sigma_{\ib,k})$ to be as before
and $\deg(\pi_11_{\ib})=0$.
{\sl In the rest of Section $1$ we'll assume $m>0$.}
Then there is a canonical inclusion of graded $\kb$-algebras
$$\cc\Rb(\Gamma)_\nu\subset\th\Rb(\Gamma)_\nu\leqno(1.2)$$
such that $1_{\ib},\varkappa_l,\sigma_k\mapsto
1_{\ib},\varkappa_l,\sigma_k$ for $\ib\in\th\! I^\nu$, $l=1,\ldots,m$,
$k=1,\ldots,m-1$ and such that $\sigma_0\mapsto\pi_1\sigma_1\pi_1$.
From now on we'll write $\sigma_0=\pi_1\sigma_1\pi_1$ whenever $m>1$.
The assignment $x\mapsto\pi_1x\pi_1$ defines 
an involution of the graded $\kb$-algebra
$\th\Rb(\Gamma)_\nu$
which normalizes 
$\cc\Rb(\Gamma)_\nu$.
Thus it yields an involution
$$\g:\cc\Rb(\Gamma)_\nu\to\cc\Rb(\Gamma)_\nu.$$
Let $\la\g\ra$ be the group of two elements generated by
$\g$. The smash product $\cc\Rb(\Gamma)_\nu\rtimes\la\g\ra$ 
is a graded $\kb$-algebra such that 
$\deg(\g)=0$. There is an unique isomorphism of graded $\kb$-algebras
$$\cc\Rb(\Gamma)_\nu\rtimes\la\g\ra\to\th\Rb(\Gamma)_\nu\leqno(1.3)$$
which is identity on $\cc\Rb(\Gamma)_\nu$ and which takes $\g$
to $\pi_1$.
%This can be checked by generators and relations.

\subhead 1.3.~The polynomial representation and the PBW
theorem\endsubhead For any $\ib$ in $\th\! I^\nu$ let $\th\Fb_\ib$ be the
subalgebra of $\cc\Rb(\Gamma)_{\nu}$ generated by $1_\ib$ and
$\varkappa_{\ib,l}$ with $l=1,2,\dots, m$. It is a polynomial algebra.
Let
$$\th\Fb_\nu=\bigoplus_{\ib\in\th\! I^\nu}\th\Fb_\ib.$$
The group $W_m$ acts on $\th\Fb_\nu$ via
$w(\varkappa_{\ib,l})=\varkappa_{w(\ib),w(l)}$ for any $w\in W_m$. Consider the fixed points set
$$\cc\Sb_\nu=(\th\Fb_\nu)^{\cc\! W_m}.$$
Regard $\th\Rb(\Gamma)_{\nu}$ and $\End(\th\Fb_\nu)$ as
$\th\Fb_\nu$-algebras via the left multiplication. In \cite{VV,
prop.~5.4} is given an injective graded $\th\Fb_\nu$-algebra
morphism $\th\Rb(\Gamma)_{\nu}\to\End(\th\Fb_\nu)$. It
restricts via (1.2) to an injective graded $\th\Fb_\nu$-algebra
morphism $$\cc\Rb(\Gamma)_{\nu}\to\End(\th\Fb_\nu).$$
Next, recall that $\cc W_m$ is the Weyl group of type $\D_m$ with simple
reflections $s_0,\ldots, s_{m-1}$. For each $w$ in $\cc W_m$ we
choose a reduced decomposition $\dot w$ of $w$. It has the following
form
$$w=s_{k_1}s_{k_2}\cdots s_{k_r},\quad 0\leqslant k_1,k_2,\dots,k_r\leqslant m-1.$$ We define an element $\sigma_{\dot w}$ in
$\cc\Rb(\Gamma)_{\nu}$ by 
$$\quad\sigma_{\dot w}=\sum_\ib 1_\ib\sigma_{\dot w} ,\quad
1_\ib\sigma_{\dot w}=\cases 1_\ib &\roman{if}\ r=0\vspace{2mm}
1_\ib\sigma_{k_1} \sigma_{k_2}\cdots\sigma_{k_r}
%&\roman{if}\ r>0,\,u=1,\vspace{2mm}
%\sigma_{\ib,\dot w}=\sigma_{\ib,\dot{wu}} \pi_{1}
&\roman{else},
\endcases\leqno(1.4)$$
Observe that the element $\sigma_{\dot w}$ may depend on the choice of the
reduced decomposition $\dot w$.

\proclaim{1.4.~Proposition} The $\kb$-algebra
$\cc\Rb(\Gamma)_{\nu}$ is a free (left or right)
$\th\Fb_\nu$-module with basis $\{\sigma_{\dot w};\,w\in
\cc W_m\}$. Its rank is $2^{m-1}m!$. The operator
$1_\ib\sigma_{\dot w}$ is homogeneous and its degree is independent
of the choice of the reduced decomposition $\dot w$.
\endproclaim
\noindent{\sl Proof :} The proof is the same as in \cite{VV,
prop.~5.5}. First, we filter the algebra $\cc\Rb(\Gamma)_{\nu}$ with
$1_\ib$, $\varkappa_{\ib,l}$ in degree 0 and $\sigma_{\ib,k}$ in
degree 1. The {\it Nil Hecke algebra} of type $\D_m$ is the
$\kb$-algebra $\cc\Nb\Hb_m$ generated by
$\bar\sigma_0,\bar\sigma_1,\dots,\bar\sigma_{m-1}$ with relations
$$\gathered
\bar\sigma_k\bar\sigma_{k'}=\bar\sigma_{k'}\bar\sigma_k\ \roman{if}\
1\leqs k<k'-1<m-1\ \roman{or}\ 0=k<k'\neq 2, \vspace{2mm}
\bar\sigma_{s_k(k)}\bar\sigma_{k}\bar\sigma_{s_k(k)}
=\bar\sigma_{k}\bar\sigma_{s_k(k)}\bar\sigma_{k},\quad
\bar\sigma_k^2=0.
\endgathered
$$
We can form the semidirect product $\th\Fb_\nu\rtimes\cc\Nb\Hb_m$,
which is generated by $1_\ib$, $\bar\varkappa_l$, $\bar\sigma_k$
with the  relations above and
$$\gathered
\bar\sigma_k\bar\varkappa_{l}=\bar\varkappa_{s_k(l)}\bar\sigma_k,
\quad \bar\varkappa_{l}\bar\varkappa_{l'}
=\bar\varkappa_{l'}\bar\varkappa_{l'}.
\endgathered
$$
One proves as in \cite{VV, prop.~5.5} that the map
$$\th\Fb_\nu\rtimes\cc\Nb\Hb_m\to\gr(\cc\Rb(\Gamma)_\nu),\quad
1_\ib\mapsto 1_{\ib},\quad \bar\varkappa_l\mapsto\varkappa_{l},\quad
\bar\sigma_{k}\mapsto\sigma_{k}.$$ is an isomorphism of
$\kb$-algebras.

\qed

\vskip3mm

Let $\th\Fb'_\nu=\bigoplus_\ib\th\Fb'_\ib,$ where $\th\Fb'_\ib$ is
the localization of the ring $\th\Fb_\ib$ with respect to the
multiplicative system generated by
$$\{\varkappa_{\ib,l}\pm \varkappa_{\ib,l'};\,
1\leqslant l\neq l'\leqslant m\}\cup
\{\varkappa_{\ib,l};\,l=1,2,\dots,m\}.$$

\proclaim{1.5.~Corollary} The inclusion 
$\cc\Rb(\Gamma)_{\nu}\subset\End(\th\Fb_{\nu})$ yields an isomorphism
of $\th\Fb'_\nu$-algebras 
$\th\Fb'_\nu\otimes_{\th\Fb_\nu}\cc\Rb(\Gamma)_{\nu}\to
\th\Fb'_\nu\rtimes \cc W_m,$ such that for each $\ib$ and each
$l=1,2,\dots,m$, $k=0,1,2,\dots, m-1$ we have
$$\aligned
&1_{\ib}\mapsto 1_{\ib},\vspace{2mm} 
&\varkappa_{\ib,l}\mapsto \varkappa_l1_{\ib},\vspace{2mm} 
&\sigma_{\ib,k}\mapsto \cases
(\varkappa_k-\varkappa_{s_k(k)})^{-1}(s_k-1)1_\ib & \text{if}\
i_k=i_{s_k(k)},\vspace{2mm}
(\varkappa_k-\varkappa_{s_k(k)})^{h_{i_{s_k(k)},i_k}}s_k1_\ib\hfill&
\text{if}\ i_k\neq i_{s_k(k)}.
\endcases
\endaligned\leqno(1.5)$$
\endproclaim

\noindent{\sl Proof:} Follows from \cite{VV, cor.~5.6} and
Proposition 1.4.

\qed \vskip3mm

Restricting the $\th\Fb_\nu$-action on $\cc\Rb(\Gamma)_{\nu}$ to
the $\kb$-subalgebra $\cc\Sb_\nu$ we get a structure of
graded $\cc\Sb_\nu$-algebra on $\cc\Rb(\Gamma)_{\nu}$.

\proclaim{1.6.~Proposition} (a) $\cc\Sb_\nu$ is isomorphic to the
center of $\cc\Rb(\Gamma)_{\nu}$. \vskip1mm

(b) $\cc\Rb(\Gamma)_{\nu}$ is a free graded module over
$\cc\Sb_\nu$ of rank $(2^{m-1}m!)^2$.
\endproclaim

\noindent{\sl Proof :} Part $(a)$ follows from Corollary 1.5. Part
$(b)$ follows from $(a)$ and Proposition 1.4.

\qed

\vskip2cm

\head 2.~Affine Hecke algebras of type D\endhead

\subhead 2.1.~Affine Hecke algebras of type D\endsubhead Fix $p$ in
$\kb^\times$. For any integer $m\geqslant 0$ we define the extended
affine Hecke algebra $\Hb_m$ of type $\D_m$ as follows. If 
$m>1$ then $\Hb_m$ is the $\kb$-algebra with 1 generated by
$$T_k,\quad X_l^{\pm 1},\quad
k=0,1,\dots,m-1,\quad l=1,2,\dots,m$$ satisfying the following
defining relations :

\vskip2mm
\itemitem{$(a)$} $X_lX_{l'}=X_{l'}X_l$,

\vskip2mm

\itemitem{$(b)$} $T_kT_{s_k(k)}T_k=T_{s_k(k)}T_kT_{s_k(k)}$,
$T_kT_{k'}=T_{k'}T_k$ if $1\leqs k<k'-1$ or $k=0$, $k'\neq 2$,

\vskip2mm

\itemitem{$(c)$}
$(T_k-p)(T_k+p^{-1})=0$,

\vskip2mm

\itemitem{$(d)$} $T_0X_1^{-1}T_0=X_2$, $T_kX_kT_k=X_{s_k(k)}$ 
if $k\neq 0$, $T_kX_l=X_lT_k$ if $k\neq 0,l,l-1$ or $k=0$, $l\neq 1,2$.

\vskip3mm

\noindent Finally, we set $\Hb_0=\kb\oplus\kb$ and $\Hb_1=\kb[X_1^{\pm 1}]$.

\vskip3mm

\subhead 2.2.~Remarks\endsubhead 
$(a)$ The extended affine Hecke
algebra $\Hb_m^\B$ of type $\B_m$ with parameters $p,q\in\kb^\times$
such that $q=1$ is generated by $P$, $T_k$, $X^{\pm1}_l$, $k=1,\ldots,
m-1$, $l=1,\ldots, m$ such that $T_k$, $X^{\pm1}_l$ satisfy the
relations as above and
$$\gathered P^2=1,\quad (PT_1)^2=(T_1P)^2,\quad PT_k=T_kP\
\roman{if}\ k\neq 1, \vspace{2mm} PX_1^{-1}P=X_1,\quad
PX_l=X_lP\ \roman{if}\ l\neq 1.\endgathered$$ 
%Set $T_0=PT_1P$ in $\Hb_m^\B$. 
See e.g., \cite{VV, sec.~6.1}.
There is an obvious $\kb$-algebra embedding $\Hb_m\subset\Hb_m^\B$.
Let $\g$ denote also the
involution $\Hb_m\to\Hb_m$, $a\mapsto PaP$. We have a
canonical isomorphism of $\kb$-algebras
$$\Hb_m\rtimes\la\g\ra\simeq\Hb_m^\B.$$
Compare Section 1.2.

\vskip1mm

$(b)$ 
%Fix an orthonormal basis $\chi_1,\chi_2,\dots,\chi_m$ of the
%reflection representation of $\cc W_m$ such that $s_0,s_1,\dots,s_{m-1}$
%are the reflection with respect to the simple roots
%$$a_0=\chi_2+\chi_1,\
%a_1=\chi_2-\chi_1,\
%\dots,
%\ a_{m-1}=\chi_m-\chi_{m-1}.$$
%Set $L=\bigoplus_{l=1}^m\ZZ\eps_l$. The
%weight lattice is
%$L^e=L+\ZZ({1\over 2}\sum_{l=1}^m\chi_l).$
Given a connected reductive group $G$ we call {\it affine Hecke algebra of G} the Hecke algebra of the extended affine Weyl group $W\ltimes P$, where $W$ is the Weyl group of $(G,T)$, $P$ is the group of
characters of $T$, and $T$ is a maximal torus of $G$.
Then $\Hb_m$ is the affine Hecke algebra of the group $SO(2m)$.
Let $\Hb_{m}^e$ be the affine Hecke algebra of the group $Spin(2m)$.
It is generated by $\Hb_m$ and an element $X_0$ such that
$$X_0^2=X_1X_2\dots X_m,\quad T_kX_0=X_0T_k\ \roman{if}\ k\neq 0,\quad 
T_0X_0X_1^{-1}X_2^{-1}T_0=X_0.$$
Thus $\Hb_m$ is the fixed point subset of the $\kb$-algebra automorphism
 of $\Hb_{m}^e$ taking $T_k, X_l$ to $T_k, (-1)^{\delta_{l,0}}X_l$ for all $k,l$.
Therefore, if $p$ is not a root of 1 the simple  $\Hb_m$-modules can be recovered from the Kazhdan-Lusztig classification
of the simple
$\Hb^e_{m}$-modules via Clifford theory, see e.g.,
 \cite{Re}.

\vskip3mm

\subhead 2.3.~Intertwiners and blocks of $\Hb_m$\endsubhead We
define
$$\gathered
\Ab=\kb[X_1^{\pm 1},X_2^{\pm 1},\dots,X_m^{\pm 1}],\quad
\Ab'=\Ab[\Sigma^{-1}],\quad \Hb'_m=\Ab'\otimes_\Ab\Hb_m,
\endgathered$$
where $\Sigma$ is the multiplicative set generated by
$$1-X_lX_{l'}^{\pm 1},\quad 1-p^2X_l^{\pm 1}X_{l'}^{\pm 1},\quad l\neq l'.$$ For
$k=0,\dots,m-1$ the intertwiner $\varphi_k$ is the element of $\Hb'_m$ given by
the following formulas
$$\varphi_k-1={X_k-X_{s_k(k)}\over
pX_k-p^{-1}X_{s_k(k)}}\,(T_k-p).\leqno(2.1)$$ The group $\cc W_m$
acts on $\Ab'$ as follows
$$\gathered(s_ka)(X_1,\dots,X_m)=a(X_1,\dots,X_{k+1},X_k,\dots,X_m)\
\roman{if}\ k\neq 0, \vspace{2mm}
(s_0a)(X_1,\dots,X_m)=a(X_2^{-1},X_1^{-1},\dots,X_m).\endgathered$$
There is an isomorphism of $\Ab'$-algebras
$$\Ab'\rtimes \cc W_m\to\Hb'_m,\quad
s_k\mapsto\varphi_k.$$ 
The semi-direct product group
$\ZZ\rtimes\ZZ_2=\ZZ\rtimes\{-1,1\}$ acts on $\kb^\times$ by
$(n,\eps):z\mapsto z^\eps p^{2n}$. Given a
$\ZZ\rtimes\ZZ_2$-invariant subset $I$ of $\kb^\times$ we denote by
$\Hb_m\text{-}\Modb_I$ the category of all $\Hb_m$-modules such that
the action of $X_1,X_2,\dots,X_m$  is locally finite with
eigenvalues in $I$. 
We associate to the set $I$ and to the element $p\in\kb^\times$ a quiver $\Gamma$ as follows.
The set of vertices is $I$, and there is one arrow $p^2i\to i$
whenever $i$ lies in $I$. We equip $\Gamma$ with an involution
$\theta$ such that $\theta(i)=i^{-1}$ for each vertex $i$ and such
that $\theta$ takes the arrow $p^2i\to i$ to the arrow
$i^{-1}\to p^{-2}i^{-1}$. {\it We'll assume that the set $I$ does
not contain $1$ nor $-1$ and that $p\neq 1,-1$}. Thus the involution
$\theta$ has no fixed points and no arrow may join a vertex of
$\Gamma$ to itself.

\subhead 2.4.~Remark\endsubhead We may assume that
$I=\pm\{p^n;\,n\in\ZZ_{\roman{odd}}\}$. See the discussion in
\cite{KM}. Then $\Gamma$ is of type $\A_\infty$ if $p$ has infinite
order and $\Gamma$ is of type $\A^{(1)}_r$ if $p^2$ is a primitive
$r$-th root of unity.

\vskip3mm

\subhead 2.5.~$\Hb_m$-modules versus
$\cc\Rb_{m}$-modules\endsubhead Assume that $m\geqs 1$. We define
the graded $\kb$-algebra
$$\gathered
\th\Rb_{I,m}=\bigoplus_\nu\th\Rb_{I,\nu},\quad
\th\Rb_{I,\nu}=\th\Rb(\Gamma)_\nu,\quad \cc\Rb_{I,m}=\bigoplus_\nu\cc\Rb_{I,\nu},\quad
\cc\Rb_{I,\nu}=\cc\Rb(\Gamma)_\nu,\vspace{2mm} \th\! I^m=\bigsqcup_\nu\th\!
I^\nu,\endgathered$$ where $\nu$ runs over the set of all dimension vectors in
$\th\NN I$ such that $|\nu|=2m$. When there is no risk of confusion
we abbreviate 
$$\th\Rb_\nu=\th\Rb_{I,\nu},\quad\th\Rb_{m}=\th\Rb_{I,m},\quad
\cc\Rb_\nu=\cc\Rb_{I,\nu},\quad\cc\Rb_{m}=\cc\Rb_{I,m}.$$ 
Note that $\th\Rb_\nu$ and $\th\Rb_m$ are the same as in \cite{VV, sec.~6.4},
with $\l=0$.
Note also that the $\kb$-algebra
$\cc\Rb_{m}$ may not have 1, because the set $I$ may be infinite.
We define $\cc\Rb_{m}\text{-}\Modb_0$ as the category of all
(non-graded) $\cc\Rb_{m}$-modules such that the elements
$\varkappa_1,\varkappa_2,\dots,\varkappa_m$ act locally nilpotently.
Let $\cc\Rb_{m}\text{-}\fModb_0$ and $\Hb_m\text{-}\fModb_I$ be
the full subcategories of finite dimensional modules in
$\cc\Rb_{m}\text{-}\Modb_0$ and $\Hb_m\text{-}\Modb_I$
respectively. Fix a formal series $f(\varkappa)$ in $\kb[[\varkappa]]$ such that
$f(\varkappa)=1+\varkappa$ modulo $(\varkappa^2)$.

\proclaim{2.6.~Theorem}We have an equivalence of categories
$$\cc\Rb_{m}\text{-}\Modb_0\to\Hb_m\text{-}\Modb_I,\quad M\mapsto M$$
which is given by \vskip1mm
\itemitem{$(a)$}
$X_l$ acts on $1_\ib M$ by $i_l^{-1}f(\varkappa_l)$ for each
$l=1,2,\dots,m$, \vskip1mm
\itemitem{$(b)$}if $m>1$ then $T_k$ acts on $1_\ib M$
as follows for each $k=0,1,\dots, m-1,$
$$\matrix
&\displaystyle
{{(pf(\varkappa_k)-p^{-1}f(\varkappa_{s_k(k)}))(\varkappa_k-\varkappa_{s_k(k)})
\over f(\varkappa_k)-f(\varkappa_{s_k(k)})}\sigma_k+p} \hfill&
\text{if}\ i_{s_k(k)}=i_k,\hfill\vspace{2mm} &\displaystyle{
{f(\varkappa_k)-f(\varkappa_{s_k(k)})\over
(p^{-1}f(\varkappa_k)-pf(\varkappa_{s_k(k)}))
({\varkappa_k}-{\varkappa_{s_k(k)})}}\sigma_k+{(p^{-2}-1)f(\varkappa_{s_k(k)})
\over pf(\varkappa_k)-p^{-1}f(\varkappa_{s_k(k)})}} &\text{if}\
i_{s_k(k)}=p^2i_{k},\hfill\vspace{2mm} &\displaystyle{
{pi_{k}f(\varkappa_k)-p^{-1}i_{s_k(k)}f(\varkappa_{s_k(k)})\over
i_{k}f(\varkappa_k)-i_{s_k(k)}f(\varkappa_{s_k(k)})}\sigma_k+
{(p^{-1}-p)i_{k}f(\varkappa_{s_k(k)})\over
i_{s_k(k)}f(\varkappa_k)-i_{k}f(\varkappa_{s_k(k)})}}
\hfill&\text{if}\  i_{s_k(k)}\neq i_{k},p^2i_{k} .\hfill
\endmatrix$$
\endproclaim

\noindent{\sl Proof :} This follows from \cite{VV, thm.~6.5} by
Section 1.2 and Remark 2.2$(a)$. One can also prove it by reproducing the
arguments in loc.~cit.~ by using (1.5) and (2.1).

\qed

\vskip3mm

\proclaim{2.7.~Corollary} There is an equivalence of categories
$$\Psi:\cc\Rb_{m}\text{-}\fModb_0\to\Hb_m\text{-}\fModb_I,\quad
M\mapsto M.$$
\endproclaim

\vskip3mm

\subhead 2.8.~Remarks\endsubhead
$(a)$
Let $\gen$ be the Lie algebra of $G=SO(2m)$. Fix a maximal torus $T\subset G$.
The group of characters of $T$ is the lattice $\bigoplus_{l=1}^m\ZZ\,\eps_l$, with Bourbaki's notation. 
Fix a dimension vector $\nu\in\th\NN I$. Recall the sequence $\ib_e=(i_{1-m},\dots,i_{m-1},i_m)$ 
from Section 0.3. 
Let $g\in T$ be the element such that $\eps_l(g)=i_l^{-1}$ for each $l=1,2,\dots,m$. Recall also the notation
$\th\Vcb_\nu$, $\Vb$, $\th\!E_\Vb$, and $\th\!G_\Vb$ from \cite{VV}. Then $\Vb$ is an object of $\th\Vcb_\nu$,
$\th\!G_\Vb=G_g$ is the centralizer of $g$ in $G$,  and
$$\th\!E_\Vb=\gen_{g,p},\quad\gen_{g,p}=\{x\in\gen;,\ad_g(x)=p^2x\}.$$
Let $F_g$ be the set of all Borel Lie subalgebras of $\gen$ 
fixed by the adjoint action of $g$.
It is a non connected manifold whose connected components are labelled by $\th\!I^\nu_+$.
In Section 3.14 we'll introduce two central idempotents
$1_{\nu, +}$, $1_{\nu, -}$ of $\cc\Rb_\nu$. This yields a graded $\kb$-algebra
decomposition
$$\cc\Rb_\nu =\cc\Rb_\nu 1_{\nu, +}\oplus\cc\Rb_\nu 1_{\nu, -}.$$
By \cite{VV, thm.~5.8} the graded $\kb$-algebra  $\cc\Rb_\nu 1_{\nu, +}$ is isomorphic to 
$$\Ext^*_{G_g}(\Lc_{g,p},\Lc_{g,p}),$$
where $\Lc_{g,p}$ is the direct image of the constant perverse sheaf by the projection
$$\{ (\ben,x)\in F_g\times\gen_{g,p};\,x\in\ben\}\to\gen_{g,p},\quad (\ben,x)\mapsto x.$$
The complex $\Lc_{g,p}$ has been extensively studied by Lusztig, see e.g., 
\cite{L1}, \cite{L2}. We hope to come back to this elsewhere.

\vskip1mm

$(b)$ The results in Section 2.5 hold true if $\kb$ is any field. Set $f(\varkappa)=1+\varkappa$ for instance.

\subhead 2.9.~Induction and restriction of
$\Hb_m$-modules\endsubhead For $i\in I$ we define  functors
$$\gathered
E_i:\Hb_{m+1}\text{-}\fModb_I\to\Hb_{m}\text{-}\fModb_I,\vspace{2mm}
F_i:\Hb_m\text{-}\fModb_I\to\Hb_{m+1}\text{-}\fModb_I,
\endgathered\leqno(2.2)$$
where $E_iM\subset M$ is the generalized $i^{-1}$-eigenspace  of
the $X_{m+1}$-action, and where
$$F_iM=\Ind_{\Hb_m\otimes\kb[X_{m+1}^{\pm 1}]}^{\Hb_{m+1}}
(M\otimes\kb_{i}).$$ Here $\kb_{i}$ is the 1-dimensional
representation of $\kb[X_{m+1}^{\pm 1}]$ defined by $X_{m+1}\mapsto
i^{-1}$.

\vskip2cm

\head 3.~Global bases of $\cc\Vb$ and projective graded
$\cc\Rb$-modules\endhead

\subhead 3.1.~The Grothendieck groups of $\upd\Rb_{m}$\endsubhead
The graded $\kb$-algebra $\upd\Rb_{m}$ is free of finite rank over
its center by Proposition 1.6, a commutative graded $\kb$-subalgebra. Therefore any
simple object of $\upd\Rb_{m}$-$\modb$ is finite-dimensional and
there is a finite number of isomorphism classes of simple modules in
$\upd\Rb_{m}$-$\modb$. The Abelian group $G(\upd\Rb_{m})$ is free
with a basis formed by the classes of the simple objects of
$\upd\Rb_{m}$-$\modb$. The Abelian group $K(\upd\Rb_{m})$ is free
with a basis formed by the classes of the indecomposable projective
objects. Both $G(\upd\Rb_{m})$ and $K(\upd\Rb_{m})$ are free $\Ac$-modules,
where $v$ shifts the grading by $1$. We consider the following $\Ac$-modules
$$\gathered
\upd\Kb_I=\bigoplus_{m\geqslant 0}\upd\Kb_{I,m},
\quad\upd\Kb_{I,m}=K(\upd\Rb_{m}),\vspace{1mm}
\upd\Gb_I=\bigoplus_{m\geqslant 0}\upd\Gb_{I,m},\quad
\upd\Gb_{I,m}=G(\upd\Rb_{m}).
\endgathered$$
We'll also abbreviate
$$\upd\Kb_{I,*}=\bigoplus_{m> 0}\upd\Kb_{I,m},\quad \upd\Gb_{I,*}=\bigoplus_{m>
0}\upd\Gb_{I,m}.$$ From now on, to unburden the notation we may
abbreviate $\upd\Rb=\upd\Rb_m$, hoping it will not create any
confusion. For any $M,N$ in $\upd\Rb$-$\modb$ we set
$$(M:N)
=\gdim(M^\omega\otimes_{\upd\Rb}N),\quad \la M:N\ra=
\gdim\,\gHom_{\upd\Rb}(M,N),$$ 
where $\omega$ is the involution defined in Section 1.1. The Cartan pairing is the perfect
$\Ac$-bilinear form 
$$\upd\Kb_I\times \upd\Gb_I\to \Ac,\quad (P,M)\mapsto\la P: M\ra.$$

First, we concentrate on the $\Ac$-module $\upd\Gb_I$. Consider
the duality
$$\upd\Rb\text{-}\fmodb\to\upd\Rb\text{-}\fmodb,\quad
M\mapsto M^\flat=\gHom(M,\kb),$$ with the action and the grading
given by
$$(xf)(m)=f(\omega(x)m),\quad(M^\flat)_d=\Hom(M_{-d},\kb).$$
This duality functor yields an $\Ac$-antilinear map
$$\upd\Gb_I\to\upd\Gb_I,\quad M\mapsto M^\flat.$$
Let $\upd\!B$ denote the set of isomorphism classes of simple objects 
of $\upd\Rb$-$\fModb_0$. We can now define the
upper global basis of $\upd\Gb_I$ as follows. The proof is given
in Section 3.21.

\proclaim{3.2.~Proposition/Definition} For each $b$ in $\upd\!B$ there
is a unique selfdual irreducible graded $\upd\Rb$-module
$\upd\!G^\up(b)$ which is isomorphic to $b$ as a (non graded)
$\upd\Rb$-module. We set $\upd\!G^\up(0)=0$ and
$\upd\Gb^\up=\{\upd\!G^\up(b);\,b\in \upd\!B\}$. Hence $\upd\Gb^\up$
is a $\Ac$-basis of $\upd\Gb_I$.
\endproclaim

Now, we concentrate on the $\Ac$-module $\upd\Kb_I$. We equip
$\upd\Kb_I$ with the symmetric $\Ac$-bilinear form
$$\gathered
\upd\Kb_I\times \upd\Kb_I\to \Ac,\quad (M,N)\mapsto(M:N).
\endgathered\leqno(3.1)$$
Consider the duality
$$\upd\Rb\text{-}\proj\to\upd\Rb\text{-}\proj,\quad
P\mapsto P^\sharp=\gHom_{\upd\Rb}(P,\upd\Rb),$$ with the action and
the grading given by
$$(xf)(p)=f(p)\omega(x),\quad(P^\sharp)_d=\Hom_{\upd\Rb}(P[-d],\upd\Rb).$$
This duality functor yields an $\Ac$-antilinear map
$$\upd\Kb_I\to\upd\Kb_I,\quad P\mapsto P^\sharp.$$
Set $\Kc=\QQ(v)$. Let $\Kc\to\Kc$, $f\mapsto \bar f$ be the unique
involution such that $\bar v=v^{-1}$.

\proclaim{3.3.~Definition} For each $b$ in $\upd\!B$ let
$\upd\!G^\low(b)$ be the unique indecomposable graded module in
$\upd\Rb$-$\proj$ whose top is isomorphic to $\upd\!G^\up(b)$. We
set $\upd\!G^\low(0)=0$ and $\upd\Gb^\low=\{\upd\!G^\low(b);\,b\in
\upd\!B\}$, a $\Ac$-basis of $\upd\!\Kb_I$.
\endproclaim

\proclaim{3.4.~Proposition} (a) We have $\la
\upd\!G^\low(b):\upd\!G^\up(b')\ra=\delta_{b,b'}$ for each $b,b'$ in
$\upd\!B$.

\vskip1mm

(b) We have $\la P^\sharp:M\ra=\overline{\la P:M^\flat\ra}$ for each
$P$, $M$.

\vskip1mm

(c) We have ${\upd\!G^\low(b)}^\sharp=\upd\!G^\low(b)$ for each $b$
in $\upd\!B$.
\endproclaim

\noindent 
The proof is the same as in \cite{VV, prop.~8.4}.

\vskip3mm

\subhead 3.5.~Example\endsubhead Set $\nu=i+\theta(i)$ and
$\ib=i\theta(i)$. Consider the graded $\upd\Rb_\nu$-modules
$$\upd\Rb_\ib=\upd\Rb 1_{\ib}=1_{\ib}\upd\Rb, \quad \upd\Lb_\ib=\top(\upd\Rb_\ib).$$
The global bases are given by
$$\upd\Gb^\low_\nu=\{\upd\Rb_\ib,\, \upd\Rb_{\theta(\ib)}\},\quad \upd\!
\Gb^\up_\nu=\{\upd\Lb_\ib,\, \upd\Lb_{\theta(\ib)}\}.$$ For $m=0$ we
have $\upd\Rb_0=\kb\oplus\kb$. Set
$\phi_+=\kb\oplus 0$ and $\phi_-=0\oplus\kb.$
We have
$$\upd\Gb^\low_0=\upd\Gb^\up_0=\{\phi_+,\,\phi_-\}.$$

\vskip3mm

\subhead 3.6.~Definition of the operators
$e_i,f_i,e'_i,f'_i$\endsubhead In this section we'll always assume
$m>0$ unless specified otherwise. First, let us introduce the
following notation. Let $D_{m,1}$ be the set of
minimal representative in $\cc W_{m+1}$ of the cosets in
$\cc W_m\backslash \cc W_{m+1}$. Write 
$$D_{m,1;m,1}=D_{m,1}\cap(D_{m,1})^{-1}.$$
For each element $w$ of $D_{m,1;m,1}$ we set
$$W(w)=\cc W_m\cap w(\cc W_m)w^{-1}.$$
Let $\Rb_1$ be the $\kb$-algebra generated by elements $1_i$,
$\varkappa_{i}$, $i\in I$, satisfying the defining relations
$1_i\,1_{i'}=\delta_{i,i'}1_i$ and
$\varkappa_{i}=1_{i}\varkappa_{i}1_i$. We equip $\Rb_1$ with the grading
such that $\deg(1_i)=0$ and $\deg(\varkappa_{i})=2$. Let
$$\Rb_i=1_i\Rb_1=\Rb_11_i,\quad\Lb_i=\top(\Rb_i)=\Rb_i/(\varkappa_i).$$
Then $\Rb_i$ is a graded projective $\Rb_1$-module and $\Lb_i$ is
simple.
We abbreviate
$$\th\Rb_{m,1}=\th\Rb_{m}\otimes\Rb_{1},\quad \upd\Rb_{m,1}=\upd\Rb_{m}\otimes\Rb_{1}.$$
There is an unique inclusion of graded $\k$-algebras
$$\gathered
\th\Rb_{m,1}\to\th\Rb_{m+1},\vspace{2mm} 1_\ib\otimes
1_{i}\mapsto1_{\ib'}, \vspace{2mm}
1_{\ib}\otimes\varkappa_{i,l}\mapsto \varkappa_{\ib',m+l},
\vspace{2mm} \varkappa_{\ib,l}\otimes 1_{i}\mapsto
\varkappa_{\ib',l}, \vspace{2mm} \pi_{\ib,1}\otimes 1_{i}\mapsto
\pi_{\ib',1}, \vspace{2mm} \sigma_{\ib,k}\otimes 1_{i}\mapsto
\sigma_{\ib',k},
\endgathered\leqno(3.2)$$
where, given $\ib\in\th\! I^m$ and $i\in I$, we have set
$\ib'=\theta(i)\ib i$, a sequence in $\th\!I^{m+1}$. This inclusion
restricts to an inclusion
$\upd\Rb_{m,1}\subset\upd\Rb_{m+1}.$

\proclaim{3.7.~Lemma} The graded $\upd\Rb_{m,1}$-module
$\upd\Rb_{m+1}$ is free of rank $2(m+1)$.
\endproclaim

\noindent{\sl Proof :} For each $w$ in $D_{m,1}$ we have the element
$\sigma_{\dot w}$ in $\upd\Rb_{m+1}$ defined in (1.5). Using
filtered/graded arguments it is easy to see that
$$\upd\Rb_{m+1}=\bigoplus_{w\in D_{m,1}}\upd\Rb_{m,1}\sigma_{\dot w}.$$

\qed

\vskip3mm

We define a triple of adjoint functors $(\psi_!,\psi^*,\psi_*)$
where
$$\psi^*\,:\,
\upd\Rb_{m+1}\text{-}\modb\to
\upd\Rb_{m}\text{-}\modb\times\Rb_{1}\text{-}\modb$$ is the
restriction and $\psi_!$, $\psi_*$ are given by
$$\aligned
&\psi_!\,:\,\cases
\upd\Rb_{m}\text{-}\modb\times\Rb_{1}\text{-}\modb\to
\upd\Rb_{m+1}\text{-}\modb,\vspace{2mm} (M,M')\mapsto
\upd\Rb_{m+1}\otimes_{\upd\Rb_{m,1}} (M\otimes M'),\endcases
\vspace{2mm} &\psi_*\,:\,\cases
\upd\Rb_{m}\text{-}\modb\times\Rb_{1}\text{-}\modb\to
\upd\Rb_{m+1}\text{-}\modb,\vspace{2mm} (M,M')\mapsto
\gHom_{\upd\Rb_{m,1}}( \upd\Rb_{m+1},M\otimes M').\endcases
\endgathered
$$
 First, note that the functors $\psi_!$, $\psi^*$, $\psi_*$ commute with the shift of the grading.
 Next, the functor $\psi^*$ is exact, and it takes finite dimensional graded modules to finite
dimensional ones. 
The right graded
$\upd\Rb_{m,1}$-module $\upd\Rb_{m+1}$ is free of finite rank.
Thus $\psi_!$ is exact, and it takes finite dimensional graded modules to finite
dimensional ones.
The left graded
$\upd\Rb_{m,1}$-module $\upd\Rb_{m+1}$ is also free of finite rank.
Thus the functor $\psi_*$ is exact, and it takes finite dimensional graded modules to finite
dimensional ones. 
Further $\psi_!$ and $\psi^*$
take projective graded modules to projective ones, because they are left adjoint to the exact functors
$\psi^*$, $\psi_*$ respectively. 
To summarize, the functors $\psi_!$, $\psi^*$, $\psi_*$ are exact and take finite dimensional graded modules to finite dimensional ones, and the functors $\psi_!$, $\psi^*$ take projective graded modules to projective ones.

For any graded $\upd\Rb_m$-module $M$we write
$$\gathered
f_i(M)
=\upd\Rb_{m+1}1_{m,i} \otimes_{\upd\Rb_{m}}M,\vspace{2mm}
e_i(M)=\upd\Rb_{m-1}\otimes_{\upd\Rb_{m-1,1}}1_{m-1,i}
M.\endgathered\leqno(3.3)$$ 
Let us explain these formulas. The symbols $1_{m,i}$ and $1_{m-1,i}$ are given by
$$1_{m-1,i}M=\bigoplus_{\ib}1_{\theta(i)\ib i}M,\quad \ib\in\th\! I^{m-1}.$$
Note that $f_i(M)$ is a graded $\upd\Rb_{m+1}$-module, while $e_i(M)$
is a graded $\upd\Rb_{m-1}$-module.
The tensor product in the definition of $e_i(M)$ is relative to the graded
$\kb$-algebra homomorphism
$$\upd\Rb_{m-1,1}\to\upd\Rb_{m-1}\otimes\Rb_1\to\upd\Rb_{m-1}\otimes\Rb_i\to
\upd\Rb_{m-1}\otimes\Lb_i=\upd\Rb_{m-1}.$$
In other words, let $e'_i(M)$ be the graded $\upd\Rb_{m-1}$-module obtained by taking the direct summand $1_{m-1,i}M$ and restricting it to $\upd\Rb_{m-1}$. Observe that if $M$ is finitely generated then
$e'_i(M)$ may not lie in $\upd\Rb_{m-1}$-$\modb$. To remedy this, since $e'_i(M)$ affords a
$\upd\Rb_{m-1}\otimes\Rb_i$-action we consider the graded $\upd\Rb_{m-1}$-module
$$e_i(M)=e'_i(M)/\varkappa_ie'_i(M).$$

\proclaim{3.8.~Definition} The functors $e_i$, $f_i$ preserve the category $\upd\Rb$-$\projb$,
yielding $\Ac$-linear operators on $\upd\Kb_{I}$ which act on $\upd\Kb_{I,*}$ by the formula $(3.3)$ and
satisfy also $$\gather f_i(\phi_+)=\upd \Rb_{\theta(i)i},\quad
f_i(\phi_-)=\upd \Rb_{i\theta(i)},\quad
e_i(\Rb_{\theta(j)j})=\delta_{i,j}\phi_++\delta_{i,\theta(j)}\phi_-.\endgather$$
Let $e_i$, $f_i$ denote also the $\Ac$-linear operators on $\upd\Gb_{I}$
which are the transpose of $f_i$, $e_i$ with respect to the Cartan
pairing. 
\endproclaim

Note that the symbols $e_i(M)$, $f_i(M)$ have a different meaning if $M$ is viewed as an element of
$\upd\Kb_I$ or if $M$ is viewed as an element of $\upd\Gb_I$. We hope this will not create any confusion.
The proof of the following lemma is the same as in \cite{VV,
lem.~8.9}.

\proclaim{3.9.~Lemma} (a) The operators $e_i$, $f_i$ on $\upd\Gb_I$ are given by
$$\gathered
e_i(M)=1_{m-1,i}M\quad
f_i(M)=%\psi_*(M,\Lb_i)=
\hom_{\upd\Rb_{m,1}} (\upd\Rb_{m+1},M\otimes\Lb_i),
\quad M\in\upd\Rb_{m}\text{-}\fmodb.
\endgathered$$

(b) For each $M\in\upd\Rb_m\text{-}\modb$,
$M'\in\upd\Rb_{m+1}\text{-}\modb$ we have
$$(e'_i(M'):M)=(M':f_i(M)).$$

\vskip1mm

(c) We have $f_i(P)^\sharp=f_i(P^\sharp)$ for each
$P\in\upd\Rb\text{-}\proj$.

\vskip2mm

(d) We have $e_i(M)^\flat=e_i(M^\flat)$ for each
 $M\in\upd\Rb\text{-}\fmodb$.
\endproclaim

\subhead 3.10.~Induction of $\Hb_m$-modules versus induction of
$\upd\Rb_{m}$-modules\endsubhead Recall the functors $E_i$, $F_i$ on
$\Hb\text{-}\fModb_I$ defined in (2.2). We have also the functors
$$\forb: \upd\Rb_m\text{-}\fmodb\to \upd\Rb_m\text{-}\fModb_0,
\quad
\Psi:\upd\Rb_{m}\text{-}\fModb_0\to\Hb_m\text{-}\fModb_I,$$
where $\forb$ is the forgetting of the grading. Finally we define
functors
$$\gathered
E_i:\upd\Rb_{m}\text{-}\fModb_0\to\upd\Rb_{m-1}\text{-}\fModb_0,\quad
E_iM=1_{m-1,i}M, \vspace{2mm}
F_i:\upd\Rb_{m}\text{-}\fModb_0\to\upd\Rb_{m+1}\text{-}\fModb_0,\quad
F_iM=\psi_!(M,\Lb_i).
\endgathered\leqno(3.4)$$

\proclaim{3.11.~Proposition} There are canonical isomorphisms of
functors
$$E_i\circ\Psi=\Psi\circ E_{i},\quad
F_i\circ\Psi=\Psi\circ F_{i},\quad E_i\circ\forb=\forb\circ
e_i,\quad F_i\circ\forb=\forb\circ f_{\theta(i)}.$$
\endproclaim

\noindent{\sl Proof :} The proof is the same as in \cite{VV, prop.~8.17}.

\qed

 \vskip3mm

\proclaim{3.12.~Proposition} (a) The functor $\Psi$ yields an
isomorphism of Abelian groups
$$\bigoplus_{m\geqslant 0}[\upd\Rb_{m}\text{-}\fModb_0]=
\bigoplus_{m\geqslant 0}[\Hb_m\text{-}\fModb_I].$$ The functors
$E_i$, $F_i$ yield endomorphisms of both sides which are intertwined
by $\Psi$.

\vskip1mm

(b) The functor $\forb$ factors to a group isomorphism
$$\upd\Gb_I/(v-1)=
\bigoplus_{m\geqslant 0}[\upd\Rb_{m}\text{-}\fModb_0].$$
\endproclaim

\noindent{\sl Proof :} Claim $(a)$ follows from Corollary 2.7 and
Proposition 3.11. Claim $(b)$ follows from Proposition 3.2.

\qed

\vskip3mm

\subhead 3.13.~ Type D versus type B\endsubhead We can compare 
the previous constructions with their analogues in 
type B. Let $\th\Kb$, $\th
\!B$, $\th\!G^\low$, etc, denote the type B analogues of
$\upd\Kb$, $\upd\!B$, $\upd\!G^\low$, etc. See \cite{VV} for details.
We'll use the same
notation for the functors $\psi^\ast$, $\psi_!$, $\psi_\ast$,
$e_i$, $f_i$, etc, on the type $\B$ side and on the type $\D$ side.
Fix $m>0$ and $\nu\in\th\NN I$ such that $|\nu|=2m$. 
The inclusion of graded $\kb$-algebras
$\upd\Rb_\nu\subset\th\Rb_\nu$ in (1.2) yields a
restriction functor $$\res:
\th\Rb_\nu\text{-}\modb\to\upd\Rb_\nu\text{-}\modb$$ and an
induction functor
$$\indu:\upd\Rb_\nu\text{-}\modb\to\th\Rb_\nu\text{-}\modb,\quad 
M\mapsto \th\Rb_\nu\otimes_{\upd\Rb_\nu}M.$$
Both functors are exact, they map finite dimensional graded modules to
finite dimensional ones, and they map projective graded modules to projective ones. Thus,
they yield morphisms of $\Ac$-modules
$$\gathered\res: \th\Kb_{I,m}\to\upd\Kb_{I,m},\quad \res:
\th\!\Gb_{I,m}\to\upd\Gb_{I,m},\vspace{2mm} \indu:
\upd\Kb_{I,m}\to\th\Kb_{I,m},\quad \indu:
\upd\Gb_{I,m}\to\th\!\Gb_{I,m}.\endgathered$$ Moreover, for any
$P\in\th\Kb_{I,m}$ and any $L\in\th\!\Gb_{I,m}$ we have
$$\gathered \res(P^\sharp)=(\res P)^\sharp,\quad
\indu(P^\sharp)=(\indu P)^\sharp\vspace{2mm} \res(L^\flat)=(\res
L)^\flat,\quad \indu(L^\flat)=(\indu L)^\flat.
\endgathered\leqno(3.5)$$
Note also that ind and res are left and right adjoint functors, because
$$\th\Rb_\nu\otimes_{\upd\Rb_\nu}M=\hom_{\upd\Rb_\nu}(\th\Rb_\nu,M)$$
as graded $\th\Rb_\nu$-modules.

\proclaim {3.14.~Definition} For any graded $\upd\Rb_\nu$-module
$M$ we define the graded $\upd\Rb_\nu$-module $M^{\g}$ with the
same underlying graded $\kb$-vector space as $M$ such that the action of
$\upd\Rb_\nu$ is twisted by $\g$, i.e.,
the graded $\kb$-algebra $\upd\Rb_\nu$ acts on $M^{\g}$ by 
$a\, m=\g(a)m$ for  $a\in\upd\Rb_\nu$ and $m\in M$. Note that
$(M^{\g})^{\g}=M$, and that $M^\g$ is
simple (resp. projective, indecomposable) if $M$ has the
same property.
\endproclaim

For any graded $\upd\Rb_m$-module $M$  we have canonical isomorphisms of
$\upd\Rb$-modules
$$(f_i(M))^{\g}= f_i(M^{\g}),\quad (e_i(M))^{\g}=
e_i(M^{\g}).$$ The first isomorphism is given
by
$$\upd\Rb_{m+1}1_{m,i}\otimes_{\upd\Rb_{m}}M\to\upd\Rb_{m+1}1_{m,i}\otimes_{\upd\Rb_{m}}M,\quad
a\otimes m\mapsto \g(a)\otimes m.$$ The second one is the
identity map on the vector space $1_{m,i}M$.

Recall that $\th\! I^\nu$ is the disjoint union of
$\th\! I^\nu_+$ and $\th\! I^\nu_-$. We set
$$1_{\nu,+}=\sum_{\ib\in\th\! I^\nu_+}1_{\ib},\quad 1_{\nu,-}=\sum_{\ib\in\th
\! I^\nu_-}1_{\ib}.$$ 

\vskip2mm

\proclaim{3.15.~Lemma} Let $M$ be a graded $\upd\Rb_\nu$-module.

\vskip2mm

\itemitem{$(a)$} Both $1_{\nu,+}$ and $1_{\nu,-}$ are central idempotents in
$\upd\Rb_\nu$. We have $1_{\nu,+}=\g(1_{\nu,-})$.

\vskip2mm
\itemitem{$(b)$} There is a decomposition of graded $\upd\Rb_\nu$-modules
$M=1_{\nu,+}M\oplus1_{\nu,-}M.$

\vskip2mm

\itemitem{$(c)$} We have a canonical
isomorphism of $\upd\Rb_\nu$-modules
$\res\circ\indu(M)=M\oplus M^{\g}.$

\vskip2mm

\itemitem{$(d)$}
If there exists $a\in\{+,-\}$ such that $1_{\nu,-a}M=0$,
then there are canonical isomorphisms of graded $\upd\Rb_\nu$-modules
$$M=1_{\nu,a}M,
\quad 
0=1_{\nu,a}M^\g,
\quad 
M^{\g}=1_{\nu,-a}M^\g.$$
\endproclaim

\noindent{\sl Proof: }Part $(a)$ follows from Proposition 1.6 and the
equality $\eps_1(\th\! I^\nu_+)=\th\! I^\nu_-$. Part $(b)$
follows from $(a)$, $(c)$ is given by definition, and $(d)$ follows from $(a)$,
$(b)$.

\qed

\vskip3mm

Now, let $m$ and $\nu$ be as before. 
Given $i\in I$, we set $\nu'=\nu+i+\theta(i)$. 
There is an obvious inclusion $W_m\subset W_{m+1}$. 
Thus the group $W_m$ acts on $\th\! I^{\nu'}$, and the map
$$\th\! I^\nu\to\th\! I^{\nu'},\quad \ib\mapsto \theta(i)\ib i$$ is
$W_m$-equivariant. Thus there is $a_i\in\{+,-\}$
such that the image of $\th\! I^\nu_+$ is contained in $\th
I^{\nu'}_{a_i}$, and the image of $\th\! I^\nu_-$ is contained in $\th
I^{\nu'}_{-a_i}$.

\proclaim{3.16.~Lemma} Let $M$ be a graded $\upd\Rb_\nu$-module such
that $1_{\nu,-a}M=0,$ with $a=+,-$. Then we have
$$1_{\nu',-a_ia}f_i(M)=0,\quad 1_{\nu',a_ia}f_{\theta(i)}(M)=0.$$
\endproclaim
\noindent{\sl Proof: } We have
$$\aligned 1_{\nu',-a_ia}f_i(M)
&=1_{\nu',-a_ia}\upd\Rb_{\nu'}1_{\nu,i}\otimes_{\upd\Rb_{\nu}}M
\vspace{2mm}
&=\upd\Rb_{\nu'}1_{\nu',-a_ia}1_{\nu,i}1_{\nu,a}\otimes_{\upd\Rb_{\nu}}M.\endaligned$$
Here we have identified $1_{\nu,a}$ with the image of
$(1_{\nu,a},1_i)$ via the inclusion (3.2). The definition of this
inclusion and that of $a_i$ yield that
$$1_{\nu',a_ia}1_{\nu,i}1_{\nu,a}=1_{\nu,a},\quad
1_{\nu',-a_ia}1_{\nu,i}1_{\nu,a}=0.$$ The first equality follows.
Next, note that for any $\ib\in \th\! I^\nu$, the sequences
$\theta(i)\ib i$ and $i\ib \theta(i)=\eps_{m+1}(\theta(i)\ib i)$
always belong to different $\cc W_{m+1}$-orbits. Thus, we have
$a_{\theta(i)}=-a_i$. So the second equality follows from the first.

\qed

\vskip3mm

Consider the following diagram
$$\xymatrix{\upd\Rb_{\nu}\text{-}\modb\times\Rb_{i}\text{-}\modb\ar@<0.5ex>[r]^{\qquad\psi_!}
\ar@<0.5ex>[d]^{\indu\times\Id\ }
&\upd\Rb_{\nu'}\text{-}\modb\ar@<0.5ex>[l]^{\qquad\psi^*}\ar@<0.5ex>[d]^{\indu}\\
\th\Rb_{\nu}\text{-}\modb\times\Rb_{i}\text{-}\modb\ar@<0.5ex>[r]^{\qquad\psi_!}\ar@<0.5ex>[u]^{\
\res\times\Id}
&\th\Rb_{\nu'}\text{-}\modb.\ar@<0.5ex>[l]^{\qquad\psi^*}\ar@<0.5ex>[u]^{\res}}$$

\proclaim{3.17.~Lemma} There are canonical isomorphisms of functors
$$\gathered
\indu\circ\psi_!=\psi_!\circ(\indu\times\Id),\quad
\psi^\ast\circ\indu=(\indu\times\Id)\circ\psi^\ast,\quad
\indu\circ\psi_\ast=\psi_\ast\circ(\indu\times\Id),\vspace{2mm}
\res\circ\psi_!=\psi_!\circ(\res\times\Id),\quad
\psi^\ast\circ\res=(\res\times\Id)\circ\psi^\ast,\quad
\res\circ\psi_\ast=\psi_\ast\circ(\res\times\Id).\endgathered$$
\endproclaim
\noindent{\sl Proof :} The functor $\indu$ is left and right
adjoint to $\res$. Therefore it is enough to prove the first two
isomorphisms in the first line. The isomorphism
$$\indu\circ\psi_!=\psi_!\circ(\indu\times\Id)$$
comes from the associativity of the induction.
Let us prove that
$$\psi^\ast\circ\indu=(\indu\times\Id)\circ\psi^\ast.$$ For any $M$ in
$\upd\Rb_{\nu'}$-$\modb$, the obvious inclusion
$\th\Rb_{\nu}\otimes\Rb_{i}\subset\th\Rb_{\nu'}$ yields a map
$$(\indu\times\Id)\,\psi^*(M)=(\th\Rb_{\nu}\otimes\Rb_{i})
\otimes_{\upd\Rb_{\nu}\otimes\Rb_{i}} \psi^*(M)
\to\psi^\ast(\th\Rb_{\nu'}\otimes_{\upd\Rb_{\nu}\otimes\Rb_{i}}
M).$$ Combining it with the obvious map
$$\th\Rb_{\nu'}\otimes_{\upd\Rb_{\nu}\otimes\Rb_{i}} M\to
\th\Rb_{\nu'}\otimes_{\upd\Rb_{\nu'}} M$$ we get a morphism of
$\th\Rb_{\nu}\otimes\Rb_{i}$-modules
$$(\indu\times\Id)\,\psi^*(M)
\to\psi^\ast\,\indu(M).$$ We need to
show that it is bijective. This is obvious because at the level of
vector spaces, the map above is given by 
$$M\oplus(\pi_{1,\nu}\otimes M)\to M\oplus (\pi_{1,\nu'}\otimes M),
\quad m+\pi_{1,\nu}\otimes n\mapsto m+\pi_{1,\nu'}\otimes n.$$ Here
$\pi_{1,\nu}$ and $\pi_{1,\nu'}$ denote the element
$\pi_1$ in $\th\Rb_\nu$ and $\th\Rb_{\nu'}$ respectively.

\qed

\vskip3mm

\proclaim{3.18.~Corollary}(a) The operators $e_i$, $f_i$ on
$\upd\Kb_{I,*}$ and on $\th\Kb_{I,*}$ are intertwined by the maps
$\indu$, $\res$, i.e., we have
$$\gathered e_i\circ\indu=\indu\circ e_i,\quad
f_i\circ\indu=\indu\circ f_i,\quad e_i\circ\res=\res\circ
e_i,\quad f_i\circ\res=\res\circ f_i.\endgathered$$

(b) The same result holds for the operators $e_i$, $f_i$ on
$\upd\Gb_{I,*}$ and on $\th\!\Gb_{I,*}$.
\endproclaim

\subhead 3.19\endsubhead Now, we concentrate on non graded
irreducible modules. First, let
$$\Res: \th\Rb_\nu\text{-}\Modb\to\upd\Rb_\nu\text{-}\Modb,\quad \Ind:
\upd\Rb_\nu\text{-}\Modb\to\th\Rb_\nu\text{-}\Modb$$ be the
(non graded) restriction and induction functors. We have
$$\forb\circ\res=\Res\circ\forb,\quad
\forb\circ\indu=\Ind\circ\forb.$$

\proclaim{3.20.~Lemma} Let $L$, $L'$ be irreducible
$\upd\Rb_\nu$-modules.

(a) The $\upd\Rb_\nu$-modules $L$ and $L^{\g}$ are not
isomorphic.

\vskip1mm

(b) $\Ind(L)$ is an irreducible $\th\Rb_\nu$-module, and
every irreducible $\th\Rb_\nu$-module is obtained in this way.

\vskip1mm

(c) $\Ind(L)\simeq\Ind(L')$ iff $L'\simeq L$ or $L^{\g}.$
\endproclaim

\noindent{\sl Proof:} For any 
$\th\Rb_\nu$-module $M\neq 0$, both $1_{\nu,+}M$ and $1_{\nu,-}M$ are
nonzero. Indeed, we have
$M=1_{\nu,+}M\oplus1_{\nu,-}M$, and we may suppose that $1_{\nu,+}M\neq 0$.
The automorphism $M\to M$, $m\mapsto\pi_1m$ takes $1_{\nu,+}M$ to
$1_{\nu,-}M$ by Lemma 3.15$(a)$. Hence $1_{\nu,-}M\neq 0$.

Now, to prove part $(a)$, suppose that $\phi:L\to L^\g$ is an isomorphism of
$\upd\Rb_\nu$-modules. We can  regard $\phi$ as a $\g$-antilinear
map $L\to L$. 
Since $L$ is irreducible, by Schur's lemma we may 
assume that $\phi^2=\Id_L$.
Then $L$ admits a $\th\Rb_\nu$-module structure
such that the $\upd\Rb_\nu$-action is as before and
$\pi_1$ acts as $\phi$. 
Thus, by the discussion above, neither $1_{\nu,+}L$ nor 
$1_{\nu,-}L$ is zero. This contradicts the fact that $L$ is an
irreducible $\upd\Rb_\nu$-module.

Parts $(b)$, $(c)$ follow from $(a)$ by Clifford theory, see e.g., 
\cite{RR, appendix}.

\qed

\vskip3mm

We can now prove Proposition 3.2.

\vskip1mm

\subhead 3.21.~Proof of Proposition 3.2\endsubhead Let $b\in\upd\!B$.
We may suppose that $b=1_{\nu,+}b$. By Lemma 3.20$(b)$ the module 
$\th\! b=\Ind(b)$ lies in $\th\!B$. By \cite{VV, prop.~8.2}
there exists a unique selfdual irreducible graded $\th\Rb$-module
$\th\!G^\up(\th\! b)$ which is isomorphic to $\th\! b$ as a non graded
module. Set $$\upd\!G^\up(b)=1_{\nu,+}\res(\th\!G^\up(\th\! b)).$$
By Lemma 3.15$(d)$ 
we have $\upd\!G^\up(b)=b$
as a non graded $\upd\Rb$-module, and
by (3.5) it is selfdual. This proves existence part of the proposition. 
To prove the uniqueness,
suppose that $M$ is another module with the same properties. Then
$\indu(M)$ is a selfdual graded $\th\Rb$-module which is isomorphic
to $\th\! b$ as a non graded $\th\Rb$-module. Thus we have
$\indu(M)=\th\!G^\up(\th\! b)$ by loc.~cit. 
By Lemma 3.15$(d)$ we have also
$$M=1_{\nu,+}\res(\th\!G^\up(\th\! b)).$$
So $M$ is isomorphic to $\upd\!G^\up(b)$.

\qed

\vskip3mm

\subhead 3.22.~The crystal operators on $\upd\Gb_I$ and
$\upd\!B$\endsubhead Fix a vertex $i$ in $I$. For each irreducible
graded $\upd\Rb_m$-module $M$ we define
$$\tilde e_i(M)=\soc\, (e_i(M)),\quad\tilde f_i(M)=\top\, \psi_!(M,\Lb_i),
\quad
\eps_i(M)=\max\{n\geqslant 0;\,e_i^n(M)\neq 0\}.$$

\vskip2mm

\proclaim{3.23.~Lemma} Let $M$ be an irreducible graded
$\upd\Rb$-module such that $e_i(M)$, $f_i(M)$ belong to
$\upd\Gb_{I,*}$. We have $$\indu(\tilde e_i(M))=\tilde e_i(\indu
(M)),\quad \indu(\tilde f_i(M))=\tilde f_i(\indu(M)),\quad
\eps_i(M)=\eps_i(\indu(M)).$$ In particular, $\tilde e_i(M)$ is irreducible or zero and
$\tilde f_i(M)$ is irreducible.
\endproclaim

\noindent{\sl Proof: } By Corollary 3.18 we have $\indu(
e_i(M))=e_i(\indu(M))$. 
Thus, since ind is an exact functor we have
$\indu(\tilde e_i(M))\subset e_i(\indu(M))$. 
Since $\indu$ is an additive functor, by Lemma 3.20$(b)$
we have indeed 
$$\indu(\tilde e_i(M))\subset\tilde e_i(\indu(M)).$$ Note
that $\indu(M)$ is irreducible by Lemma 3.20$(b)$, thus $\tilde
e_i(\indu(M))$ is irreducible by \cite{VV, prop.~8.21}. Since
$\indu(\tilde e_i(M))$ is nonzero, the inclusion is an isomorphism.
The fact that $\indu(\tilde e_i(M))$ is irreducible implies in
particular that $\tilde e_i(M)$ is simple. The proof of the second
isomorphism is similar. The third equality is obvious.

\qed

\vskip3mm

Similarly, for each irreducible $\upd\Rb$-module $b$ in $\upd\!B$ we define
$$\tilde E_i(b)=\soc(E_i(b)),\quad\tilde F_i(b)=\top(F_i(b)),\quad
\eps_i(b)=\max\{n\geqslant 0;\,E_i^n(b)\neq 0\}.$$ 
Hence we have 
$$\forb\circ\tilde e_i=\tilde E_i\circ\forb,\quad\forb\circ\tilde
f_i=\tilde F_i\circ\forb,\quad\eps_i=\eps_i\circ\forb.$$

\proclaim{3.24.~Proposition} For each $b,b'$ in $\upd\!B$ we have

\vskip1mm

(a) $\tilde F_i(b)\in \upd\!B$,

\vskip1mm

(b) $\tilde E_i (b)\in \upd\!B\cup\{0\}$,

\vskip1mm

(c) $\tilde F_i(b)=b'\iff\tilde E_i(b')=b$,

\vskip1mm

(d) $\eps_i(b)=\max\{n\geqslant 0;\tilde E_i^n(b)\neq 0\}$,

\vskip1mm

(e) $\eps_i(\tilde F_i(b))=\eps_i(b)+1$,

\vskip1mm

(f) if $\tilde E_i(b)=0$ for all $i$ then $b=\phi_{\pm}$.
\endproclaim

\noindent{\sl Proof:} Part $(c)$ follows from adjunction. The other
parts follow from \cite{VV, prop.~3.14} and Lemma 3.23.

\qed

\vskip3mm

\subhead 3.25.~Remark\endsubhead For any $b\in\upd\!B$ and any 
$i\neq j$ we have $\tilde F_i(b)\neq\tilde F_j(b)$. This is
obvious if $j\neq \theta(i)$. Because in this case $\tilde F_i(b)$
and $\tilde F_j(b)$ are $\upd\Rb_\nu$-modules for different $\nu$.
Now, consider the case $j=\theta(i)$. We may suppose that $\tilde
F_i(b)=1_{\nu,+}\tilde F_i(b)$ for certain $\nu$. Then by Lemma 3.16
we have $1_{\nu,+}\tilde F_{\theta(i)}(b)=0$. In particular $\tilde
F_i(b)$ is not isomorphic to $\tilde F_{\theta(i)}(b)$.

\vskip5mm

\subhead 3.26.~The algebra $\th\!\Bcb$ and the $\th\!\Bcb$-module 
$\upd\Vb$\endsubhead Following \cite{EK1,2,3} we define a
$\Kc$-algebra $\th\!\Bcb$ as follows.

\proclaim{3.27.~Definition} Let $\th\!\Bcb$ be the $\Kc$-algebra
generated by $e_i$, $f_i$ and invertible elements $t_i$, $i\in I$,
satisfying the following defining relations

\vskip2mm
\itemitem{$(a)$}
$t_it_j=t_jt_i$ and  $t_{\theta(i)}=t_i$ for all $i,j$,

\vskip2mm
\itemitem{$(b)$}
$t_ie_jt_i^{-1}=v^{i\cdot j+\theta(i)\cdot j}e_j$ and
$t_if_jt_i^{-1}= v^{-i\cdot j-\theta(i)\cdot j}f_j$ for all $i,j$,

\vskip2mm
\itemitem{$(c)$}
$e_if_j=v^{-i\cdot j}f_je_i+ \delta_{i,j}+\delta_{\theta(i),j}t_i$
for all $i,j$,

\vskip2mm
\itemitem{$(d)$} 
${\ds\sum_{a+b=1-i\cdot j}(-1)^ae_i^{(a)}e_je_i^{(b)}=
\sum_{a+b=1-i\cdot j}(-1)^af_i^{(a)}f_jf_i^{(b)}=0}$
if $i\neq j$. 
\endproclaim

\vskip2mm

\noindent Here and below we use the following notation 
$$\theta^{(a)}=\theta^a/\la a\ra!,\quad
\la a\ra=\sum_{l=1}^av^{a+1-2l},\quad
\la a\ra!=\prod_{l=1}^m\la l\ra.$$

\noindent We can now construct a representation of $\th\!\Bcb$ as
follows. By base change, the operators $e_i$, $f_i$ in Definition 3.8
yield $\Kc$-linear operators on the $\Kc$-vector space
$$\upd\Vb=\Kc\otimes_\Ac\upd\Kb_I.$$ 
We equip $\upd\Vb$ with the $\Kc$-bilinear form given by 
$$(M:N)_{_{KE}}=(1-v^2)^m\,(M:N),\quad\forall M,N\in\upd\Rb_m\text{-}\projb.$$

\proclaim{3.28.~Theorem} (a) The operators $e_i$, $f_i$ define a
representation of $\th\!\Bcb$ on $\upd\Vb$. The $\th\!\Bcb$-module
$\upd\Vb$ is generated by linearly independent vectors $\phi_+$ and
$\phi_-$ such that for each $i\in I$ we have
$$e_i\phi_\pm=0,\quad
t_i\phi_\pm=\phi_\mp,\quad \{x\in \upd\Vb;\,e_jx=0,\,\forall
j\}=\kb\, \phi_+\oplus\kb\, \phi_-.$$

(b) The symmetric bilinear form on $\upd\Vb$ is non-degenerate. We
have ${(\phi_a:\phi_{a'})_{_{KE}}=\delta_{a,a'}}$ for $a,a'=+,-,$ and
$(e_ix:y)=(x:f_iy)_{_{KE}}$ for $i\in I$ and $x,y\in\upd\Vb$.

%\vskip2mm

%(c) There is a unique $\Kc$-antilinear map $\upd\Vb\to\upd\Vb$ such
%that $P\mapsto P^\sharp$ for all graded projective module $P$. It is
%the unique $\Kc$-antilinear map such that $\phi_\pm^\sharp=\phi_\pm$
%and $(f_ix)^\sharp=f_i (x^\sharp)$ for all $x\in\upd\Vb$.

\endproclaim

\noindent{\sl Proof :} For each $i$ in $I$ we define the
$\Ac$-linear operator $t_i$ on $\upd\Kb_I$ by setting
$$t_i\phi_\pm=\phi_\mp\quad\and\quad t_i P=v^{-\nu\cdot(i+\theta(i))}P^{\g},\quad
\forall P\in\upd\Rb_{\nu}\text{-}\proj.$$ We must prove that the
operators $e_i$, $f_i$, and $t_i$ satisfy the relations of
$\th\!\Bcb$. The relations $(a)$, $(b)$ are obvious.
The relation $(d)$ is standard. It remains to check $(c)$. For
this we need a version of the Mackey's induction-restriction
theorem. Note that for $m>1$ we have
$$\gathered D_{m,1;m,1}=\{e,s_{m},\eps_{m+1}\eps_1\},\vspace{2mm}
W(e)=\cc W_m,\quad W(s_{m})=\cc W_{m-1},\quad
W(\eps_{m+1}\eps_1)=\cc W_m.\endgathered$$ Recall also that for
$m=1$ we have set $\cc W_1=\{e\}$.

\proclaim{3.29.~Lemma} Fix $i$, $j$ in $I$. Let $\mu$, $\nu$ in
$\th\NN I$ be such that $\nu+i+\theta(i)=\mu+j+\theta(j)$. Put
$|\nu|=|\mu|=2m$. The graded
$(\upd\Rb_{m,1},\upd\Rb_{m,1})$-bimodule
$1_{\nu,i}\upd\Rb_{m+1}1_{\mu,j}$ has a filtration by graded
bimodules whose associated graded is isomorphic to
$$\delta_{i,j}\bigl(\upd\Rb_\nu\otimes\Rb_i\bigr)\oplus\delta_{\theta(i),j}
\bigl((\upd\Rb_\nu)^{\g}\otimes\Rb_{\theta(i)}\bigr)[d']\oplus
A[d],$$ where $A$ is equal to
$$\matrix
&(\upd\Rb_m1_{\nu',i}\otimes \Rb_{i})\otimes_{\Rb'}
(1_{\nu',i}\upd\Rb_m\otimes \Rb_i)\hfill&\text{if}\
m>1,\hfill\vspace{2mm}
&(\upd\Rb_{\theta(\jb)}\otimes\Rb_i\otimes_{\upd\Rb_1\otimes\Rb_1}
\upd\Rb_{\theta(\ib)}\otimes\Rb_j)\oplus
(\upd\Rb_{\jb}\otimes\Rb_i\otimes_{\upd\Rb_1\otimes\Rb_1}
\upd\Rb_{\ib}\otimes\Rb_j)\hfill
&\text{if}\ m=1.\hfill\endmatrix$$ Here we have set
$\nu'=\nu-j-\theta(j)$, $\Rb'= \upd\Rb_{m-1,1}\otimes\Rb_1$,
$\ib=i\theta(i)$, $\jb=j\theta(j)$, $d=-i\cdot j$, and
${d'=-\nu\cdot(i+\theta(i))/2.}$
\endproclaim

\noindent
The proof is standard and is left to the reader.
Now, recall that for $m>1$ we have
$$f_j(P)=\upd\Rb_{m+1}1_{m,j}\otimes_{\upd\Rb_{m,1}}(P\otimes\Rb_1),\quad
e'_i(P)=1_{m-1,i}P,$$ where $1_{m-1,i}P$ is regarded as a
$\upd\Rb_{m-1}$-module. Therefore we have
$$\gathered
e'_if_j(P)=1_{m,i}\upd\Rb_{m+1}1_{m,j}\otimes_{\upd\Rb_{m,1}}(P\otimes\Rb_1),
\vspace{2mm} f_je'_i(P)=\upd\Rb_{m}1_{m-1,j}\otimes_{\upd\Rb_{m-1,1}
}(1_{m-1,i}P\otimes\Rb_1).
\endgathered$$
Therefore, up to some filtration we have the following identities

\vskip2mm

\itemitem{$\bullet$}  $e'_if_i(P)=P\otimes\Rb_i+f_ie'_i(P)[-2]$,

\vskip2mm

\itemitem{$\bullet$}
$e'_if_{\theta(i)}(P)=P^{\g}\otimes\Rb_{\theta(i)}
[-\nu\cdot(i+\theta(i))/2] +f_{\theta(i)}e'_i(P)[-i\cdot\theta(i)]$,

\vskip2mm

\itemitem{$\bullet$}  $e'_if_j(P)=f_je'_i(P)[-i\cdot j]$ if $i\neq j,\theta(j)$.
\vskip2mm

\noindent
These identities also hold for $m=1$ and $P=\upd\Rb_{\theta(i)i}$
for any $i\in I$. The first claim of part $(a)$ follows because
$\Rb_i=\kb\oplus\Rb_i[2]$. The fact that $\upd\Vb$ is
generated by $\phi_\pm$ is a corollary of Proposition~3.31 below.
Part $(b)$ of the theorem follows from \cite{KM, prop.~2.2(ii)} and
Lemma 3.9$(b)$.

\qed

\vskip3mm

\subhead 3.30.~Remarks\endsubhead 
$(a)$ The $\th\!\Bcb$-module $\cc\Vb$ is the same as the
$\th\!\Bcb$-module $V_\theta$ from \cite{KM, prop.~2.2}.
 The involution 
$\sigma:\upd\Vb\to\upd\Vb$ in \cite{KM, rem.~2.5(ii)}
is given by $\sigma(P)=P^{\g}$. It yields an involution of $\cc B$ in the obvious way.
Note that Lemma~3.20$(a)$ yields
$\sigma(b)\neq b$ for any $b\in\upd\!B$.

\vskip2mm

$(b)$
Let $\th\Vb$ be the $\th\!\Bcb$-module $\Kc\otimes_{\Ac}\th\Kb_I$
and let $\phi$ be the class of the trivial $\th\Rb_0$-module $\kb$,
see \cite{VV, thm.~8.30}. We have an inclusion of
$\th\!\Bcb$-modules
$$\th\Vb\to \upd\Vb,
\quad\phi\mapsto\phi_+\oplus\phi_-,\quad P\mapsto\res (P).$$

\proclaim{3.31.~Proposition} For any $b\in\upd\!B$ the following holds
%\itemitem{$(a)$}we have
$$\left\{\gathered
f_i(\upd\!G^\low(b))=\la\eps_i(b)+1\ra\,\upd\!G^\low(\tilde
F_ib)+\sum_{b'}f_{b,b'}\upd\!G^\low(b'),\vspace{2mm}
b'\in\upd\!B,\quad \eps_i(b')>\eps_i(b)+1,\quad f_{b,b'}\in
v^{2-\eps_i(b')}\ZZ[v],\endgathered\right.\leqno(a)$$

\vskip2mm

%\itemitem{$(b)$}we have
$$\left\{\gathered
e_i(\upd\!G^\low(b))=v^{1-\eps_i(b)}\,\upd\!G^\low(\tilde
E_ib)+\sum_{b'}e_{b,b'}\upd\!G^\low(b'),\vspace{2mm}
b'\in\upd\!B,\quad \eps_i(b')\geqs\eps_i(b),\quad e_{b,b'}\in
v^{1-\eps_i(b')}\ZZ[v].\endgathered\right.\leqno(b)$$
\endproclaim

\noindent{\sl Proof: } We prove part $(a)$, the proof for $(b)$ is
similar. If $\upd\!G^\low(b)=\phi_\pm$ this is obvious. So we
assume that $\upd\!G^\low(b)$ is a $\upd\Rb_m$-module for $m\geqs
1$. Fix $\nu\in\th\NN I$ such that $f_i(\upd\!G^\low(b))$
is a $\upd\Rb_\nu$-module. We'll abbreviate $1_{\nu,a}=1_a$ for
$a\in\{+,-\}$. Since $\upd\!G^\low(b)$ is indecomposable, it
fulfills the condition of Lemma 3.16. So there exists $a\in\{+,-\}$
such that $1_{-a}f_i(\upd\!G^\low(b))=0$. Thus, by Lemma 3.15$(c),(d)$ 
and Corollary 3.18 we have 
$$f_i(\upd\!G^\low(b))=1_a\res\,\indu
f_i(\upd\!G^\low(b))=1_a\res\,f_i\,\indu(\upd\!G^\low(b)).$$ Note that
$\th\! b=\Ind(b)$ belongs to $\th\!B$ by Lemma 3.20$(b)$. Hence (3.5) yields
$$\indu(\upd\!G^\low(b))=\th\!G^\low(\th\! b).$$ We deduce that
$$f_i(\upd\!G^\low(b))=1_a\res f_i(\th\!G^\low(\th\! b)).$$
Now, write
$$f_i(\th\!G^\low(\th\! b))=
\sum f_{\th\! b,\th\! b'}\,\th\!G^\low(\th\! b'),\quad
\th\! b'\in\th\!B.$$ Then we have 
$$f_i(\upd\!G^\low(b))=\sum f_{\th\!
b,\th\! b'}1_a\res(\th\!G^\low(\th\! b')).$$ For any $\th\!
b'\in\th\!B$ the $\upd\Rb$-module $1_a\Res(\th\! b')$ belongs to 
$\upd\!B$. Thus we have $$1_a\res(\th\!G^\low(\th\!
b'))=\upd\!G^\low(1_a\Res(\th\! b')).$$ If $\th\! b'\neq \th\! b''$
then $1_a\Res(\th\! b')\neq 1_a\Res(\th\! b'')$,
because $\th\! b'=\Ind(1_a\Res(\th\! b'))$.
Thus 
$$f_i(\upd\!G^\low(b))=
\sum f_{\th\! b,\th\! b'}\,\upd\!G^\low(1_a\Res (\th\! b')),$$
and this
is the expansion of the lhs in the lower global basis.
Finally, we have 
$$\eps_i(1_a\Res(\th\! b'))=\eps_i(\th\! b')$$ by
Lemma~3.23. So part $(a)$ follows from \cite{VV,
prop.~10.11$(b)$, 10.16}.

\qed

\vskip3mm

\subhead 3.32.~The global bases of $\upd\Vb$\endsubhead Since the
operators $e_i$, $f_i$ on $\upd\Vb$ satisfy the relations
$e_if_i=v^{-2}f_ie_i+1$, we can define the modified root operators
$\tilde{\eb}_i$, $\tilde{\fb}_i$ on the $\th\!\Bcb$-module $\upd\Vb$
as follows. For each $u$ in $\upd\Vb$ we write
$$\gathered
u=\sum_{n\geqslant 0}f_i^{(n)}u_n\ \roman{with}\ e_iu_n=0,
\vspace{1mm} \tilde{\eb}_i(u)=\sum_{n\geqslant
1}f_i^{(n-1)}u_n,\quad \tilde{\fb}_i(u)=\sum_{n\geqslant
0}f_i^{(n+1)}u_n.
\endgathered$$
Let $\Rc\subset\Kc$ be the set of functions which are regular at
$v=0$. Let $\upd\Lb$ be the $\Rc$-submodule of $\upd\Vb$ spanned by
the elements $\tilde{\fb}_{i_1}\dots\tilde{\fb}_{i_l}(\phi_\pm)$
with $l\geqslant 0$, $i_1,\dots,i_l\in I$. The following is the main
result of the paper.

\proclaim{3.33.~Theorem} (a) We have
$$\gathered
\upd\Lb=\bigoplus_{b\in \upd\!B}\Rc\, \upd\!G^\low(b),\quad
\tilde{\eb}_i(\upd\Lb)\subset\upd\Lb,
\quad\tilde{\fb}_i(\upd\Lb)\subset\upd\Lb,\vspace{2mm}
\tilde{\eb}_i(\upd\!G^\low(b))=\upd\!G^\low(\tilde E_i(b))\ \mod\
v\,\upd\Lb,\quad \tilde{\fb}_i(\upd\!G^\low(b))=\upd\!G^\low(\tilde
F_i(b))\ \mod\ v\,\upd\Lb.
\endgathered$$

(b) The assignment $b\mapsto \upd\!G^\low(b)\ \mod\ v\,\upd\Lb$
yields a bijection from $\upd\!B$ to the subset of
$\upd\Lb/v\upd\Lb$ consisting of the
$\tilde{\fb}_{i_1}\dots\tilde{\fb}_{i_l}(\phi_\pm)$'s. Further
$\upd\!G^\low(b)$ is the unique element $x\in\upd\Vb$ such that
$x^\sharp=x$ and $x=\upd\!G^\low(b)\ \mod\ v\,\upd\Lb.$

\vskip2mm

(c) For each $b,b'$ in $\upd\!B$ let $E_{i,b,b'}, F_{i,b,b'}\in\Ac$
be the coefficients of $\upd\!G^\low(b')$ in
$e_{\theta(i)}(\upd\!G^\low(b))$, $f_i(\upd\!G^\low(b))$ respectively.
Then we have
$$\gathered
E_{i,b,b'}|_{v=1}=
[F_i\Psi\forb(\upd\!G^\up(b')):\Psi\forb(\upd\!G^\up(b))],\vspace{2mm}
F_{i,b,b'}|_{v=1}=[E_i\Psi\forb(\upd\!G^\up(b')):\Psi\forb(\upd\!G^\up(b))].
\endgathered
$$
\endproclaim

\noindent{\sl Proof :} Part $(a)$ follows from \cite{EK3, thm.~4.1,
cor.~4.4}, \cite{E, Section 2.3}, and Proposition 3.31. The first claim
in $(b)$ follows from $(a)$. The second one is obvious. Part $(c)$
follows from Proposition 3.11. More precisely, by duality we can
regard $E_{i,b,b'}$, $F_{i,b,b'}$ as the coefficients of
$\upd\!G^\up(b)$ in $f_{\theta(i)}(\upd\!G^\up(b'))$ and
$e_i(\upd\!G^\up(b'))$ respectively. Therefore, by Proposition 3.11 we
can regard $E_{i,b,b'}|_{v=1}$, $F_{i,b,b'}|_{v=1}$ as the
coefficients of $\Psi\forb(\upd\!G^\up(b))$ in
$F_i\Psi\forb(\upd\!G^\up(b'))$ and $E_i\Psi\forb(\upd\!G^\up(b'))$
respectively.

\qed


\vskip2cm



 \Refs
\widestnumber\key{ABC}

\ref\key{A}\by Ariki, S. \paper On the decomposition numbers of the Hecke algebra of
$G(m,1,n)$\jour J. Math. Kyoto Univ. \vol 36\yr 1996\pages 789-808\endref

\ref\key{E}\by Enomoto, N.\paper A quiver construction of symmetric
crystals \jour Int. Math. Res. Notices\yr 2009\pages 2200-2247 \vol
12\endref

\ref\key{EK1}\by Enomoto, N., Kashiwara, M.\paper Symmetric Crystals
and the affine Hecke algebras of type B\jour Poc. Japan Acad. Ser. A Math. Sci.
\vol 83 \yr 2007\pages 135-139\endref

\ref\key{EK2}\by Enomoto, N., Kashiwara, M.\paper Symmetric Crystals
and LLT-Ariki type conjectures for the affine Hecke algebras of type
B\inbook
Combinatorial representation theory and related topics
\pages 1-20 \bookinfo RIMS K\^oky\^uroku Bessatsu, B8 
\publ Res. Inst. Math. Sci. (RIMS), Kyoto \yr 2008 \endref



\ref\key{EK3}\by Enomoto, N., Kashiwara, M.\paper Symmetric Crystals
for $\gen\len_\infty$ \jour Publications of the Research Intsitute
for Mathematical Sciences, Kyoto University\vol 44 \yr 2008\pages
837-891\endref

\ref\key{G}\by Grojnowski, I.\paper Representations of affine Hecke algebras (and affine quantum 
$GL_n$) at roots of unity\jour Internat. Math. Res. Notices\yr 1994,
\vol 5\pages 215\endref

\ref\key{KL}\by Khovanov, M., Lauda, A. D.\paper A diagrammatic
approach to categorification of quantum groups I\jour 
Represent. Theory\vol 13\yr 2009\pages 309-347\endref





\ref\key{KM}\by Kashiwara, M., Miemietz, V. \paper Crystals and
affine Hecke algebras of type $\D$\jour
Proc. Japan Acad. Ser. A Math. Sci. \vol 83  \yr 2007\pages 135-139\endref


\ref\key{L1}\by Lusztig, G.\paper Study of perverse sheaves arising from
graded Lie algebras 
\jour Advances in Math\yr
1995\pages 147-217\vol 112\endref

\ref\key{L2}\by Lusztig, G.\paper Graded Lie algebras and intersection 
cohomology
\jour  arXiv:0604.535
\endref

\ref\key{M} \by Miemietz, V.\paper
On representations of affine Hecke algebras of type $B$  
\jour Algebr. Represent. Theory  
\vol 11  
\yr 2008 
\pages 369-405
\endref





\ref\key{Re}\by Reeder, M.\paper Isogenies of Hecke algebras and a Langlands 
correspondence for ramified principal series representations
\jour Representation Theory\vol 6\yr 2002\pages 101-126\endref

\ref\key{R}\by Rouquier, R. \paper 2-Kac-Moody algebras \jour
 arXiv:0812.5023\endref

\ref\key{RR}\by Ram, A., Ramagge, J. \paper Affine Hecke algebras,
cyclotomic Hecke algebras and Clifford theory
\inbook A tribute to C. S. Seshadri (Chennai, 2002) 
\pages 428-466 
\bookinfo Trends Math. 
\publ Birkh\"auser, Basel 
\yr 2003
\endref



\ref\key{VV}\by Varagnolo, M., Vasserot, E. \paper Canonical bases
and affine Hecke algebras of type $\B$\jour arXiv:0911.5209
\endref

\endRefs
\enddocument